\documentclass[letterpaper,12pt,leqno,oneside]{article}
\usepackage{color}
\usepackage{bm}
\usepackage{exscale}
\usepackage{amsmath}
\usepackage{amsfonts}
\usepackage{stmaryrd}
\usepackage{amscd}
\usepackage{graphicx}
\usepackage{amsxtra}
\usepackage{amssymb}
\usepackage{theorem}
\usepackage[final]{epsfig}
\usepackage{eqnarray}
\usepackage{hyperref}
\newtheorem{proposition}{Proposition}[subsection]
\newtheorem{definition}[proposition]{Definition}

\newtheorem{lemma}[proposition]{Lemma}

{\theorembodyfont{\rmfamily}\newtheorem{remark}[proposition]{Remark}}
\newtheorem{theorem}[proposition]{Theorem}

{\theorembodyfont{\rmfamily}}

\newfont{\abc}{cmtt10 scaled 1200}

\def\R{\mathbb{R}}

\def\Z{\mathbb{Z}}

\def\P{\mathbb{P}}

\def\P{\mathbb{P}}

\def\TP{\mathbb{TP}}

\def\ve{\varepsilon}

\def\ra{\rightarrow}
\def\cs{\symbol{35}}
\def\p{\partial}
\def\qed{\hfill $\Box$ \\}
\def\qeda{\hfill $\Box$}
\def\mm{\mbox}
\def\v{= \emptyset}
\def\n{\neq \emptyset}

\def\M{\mathbb{A}}

\def\bp{\langle A \rangle}
\def\db{d^{\,\flat}}

\def\si{$\mathcal{S}$}

\def\bp{\langle A \rangle}

\begin{document}

\vspace*{1cm}

\begin{center}\Large{\bf{Potential Theory on Minimal Hypersurfaces I: Singularities as Martin Boundaries}}\\
\bigskip
\large{\bf{Joachim Lohkamp}\\
\medskip}

\end{center}

\noindent Mathematisches Institut, Universit\"at M\"unster, Einsteinstrasse 62, Germany\\
 {\small{\emph{e-mail: j.lohkamp@uni-muenster.de}}}

{\small {\center \tableofcontents}

\vspace{0.3cm} \setcounter{section}{1}
\renewcommand{\thesubsection}{\thesection}
%
%
%
\subsection{Introduction}\label{intro}
Area minimizing hypersurfaces and, more generally, almost minimizing hypersurfaces occur in geometry, dynamics and physics. Examples are horizons of black holes, level sets of geometric flows, and minimizing hypersurfaces used in scalar curvature geometry. A central problem is that a general (almost) minimizing hypersurface $H$ contains a complicated singularity set $\Sigma\subset H$ and $H \setminus \Sigma$ degenerates towards $\Sigma$ in a rather delicate way. Moreover many of the elliptic operators, one typically studies on $H$, also degenerate in their own way while we approach $\Sigma$.\\

In view of these entanglements the program of this paper (and its second part ~\cite{L2})  may appear surprising: we  develop a detailed potential theory on $H \setminus \Sigma$ applicable to a large class of linear elliptic second order operators. We even get a fine control over their analysis near $\Sigma$.\\

To this end, we first derive \emph{boundary Harnack inequalities} where we regard $\Sigma$ as a boundary of the open manifold $H \setminus \Sigma$. We use these central inequalities to deduce a variety of further results. For instance, we observe that $\Sigma$ is \emph{homeomorphic} to the \emph{Martin boundary}. Moreover, each boundary point is minimal. For operators \emph{naturally} associated with such a hypersurface $H$, for instance the Jacobi field operator or the conformal Laplacian, we can go even further than Martin theory on $H \setminus \Sigma$. We get \emph{stable} boundary Harnack inequalities which apply, with the same Harnack constants, to all blow-up hypersurfaces we get from infinite scalings around singular points. This considerably refines the asymptotic analysis towards $\Sigma$  by dimensional reductions we get from tangent cone approximations.\\

\noindent\textbf{\si-structures.}
This unexpected degree of control is due to a geometric property of $H \setminus \Sigma$, namely its \emph{\si-uniformity}~\cite{L1}. For the curved space $H \setminus \Sigma$ this \si-uniformity is the counterpart to \emph{uniformity} for Euclidean domains. The \si-uniformity implies the existence of a canonical \emph{conformal hyperbolic unfolding} of $H \setminus \Sigma$ into some complete \emph{Gromov hyperbolic space of  bounded geometry}, its \si-geometry. The Gromov boundary of the unfolding is homeomorphic to $\Sigma$.  For many elliptic operators, the potential theory on such hyperbolic manifolds is remarkably transparent, due to work of Ancona, cf. \cite{An1}, \cite{An2} and \cite{KL}. We apply Ancona's theory to the hyperbolic unfolding of $H \setminus \Sigma$ and use this as our starting point to study the potential theory on the original non-complete space $H \setminus \Sigma$ and the asymptotic analysis towards $\Sigma$\\

In the case of Euclidean domains one knows that, conversely, uniformity is also a \emph{necessary} prerequisite for both, the existence of hyperbolizing conformal deformations, \cite{BHK}, and due to Aikawa \cite{Ai3}, also for the validity of boundary Harnack inequalities for the Laplace operator. \\

All this makes \si-uniformity of $H\setminus\Sigma$, rather than the internal structure of $\Sigma$, the natural  input for an understanding of the asymptotic analysis near $\Sigma$. In fact, in this paper and its follow-up \cite{L2}, the structure of $\Sigma$ itself is irrelevant.\\

\noindent\textbf{Parts I and II.}
We develop this theory in two papers. In the current paper, we employ \si-structures to set up the basic potential theory for a broad class of operators, the \emph{\si-adapted operators}. In particular, we can establish surprisingly robust boundary Harnack inequalities for these. We apply them to derive a Martin theory and solve classical boundary value problems for the boundary $\Sigma$. We also consider eigenvalue problems when the operator is symmetric. For this purpose we introduce \si-Sobolev spaces and other variational tools.\\

In the sequel~\cite{L2} we shall unravel another intrinsic property of almost minimizers. Namely, they carry a \emph{Hardy structure}. In simple terms, Hardy structures provide us with a handy tool to prove the \si-adaptedness of various classical operators on area minimizers. Furthermore, we show that for a naturally defined \si-adapted operator $L$, \emph{minimal growth} of positive solutions of $L\phi=0$ towards $\Sigma$ is a \emph{stable} property. It persists under perturbations or blow-ups of the underlying spaces. With these results we can develop a dimensional induction scheme for the asymptotic analysis of \si-adapted operators near $\Sigma$.\\

\noindent\textbf{Remark.} This work is a substantially extended version of (parts of) earlier lecture notes by the author \cite{L} available through the arxiv. These notes may be used as a panorama of how this potential theory is linked to other geometro-analytic topics related to minimal surfaces.

%
\subsubsection{Basic Concepts and Notations}\label{notation}
In this paper $H^n$ denotes a connected integer multiplicity rectifiable current of dimension $n \ge 2$  inside some complete smooth Riemannian manifold $(M^{n+1}, g_M)$. (For some facts from geometric measure theory which we will use in the sequel we refer to \cite[Appendix A]{L1}.)\\

We denote the set of singular points of $H$ by $\Sigma_H$ or simply $\Sigma$ if there is no risk of confusion. In the special case of a minimal cone $C$ we also write $\sigma_C$ instead of $\Sigma_C$ whenever we view $C$ as a tangent cone. The induced Riemannian metric on $H\setminus\Sigma$ will be denoted by $g_H$. Furthermore, we refer to the distance function $d_H$ on $H$ inherited from the ambient space as the \emph{intrinsic distance}. Since $H$ degenerates towards $\Sigma$, it is a non-obvious fact that $d_H$ is a distance function on $H$ and that, for compact $H$, the diameter is finite, cf.~\cite[Theorem 1.8]{L1} for details.\\

For $\lambda>0$ we let $\lambda\cdot M$ 
denote the conformally rescaled Riemannian manifolds $(M,\lambda^2 \cdot  g)$; we write $\lambda\cdot H$ for the corresponding hypersurface.\\

\noindent\textbf{Area Minimizers.} \,
We shall consider the following sets of connected integer multiplicity rectifiable currents $H^n$.
\begin{description}
  \item[${\cal{H}}^c_n$:] $H^n \subset M^{n+1}$ is  compact locally mass minimizing without boundary.
  \item[${\cal{H}}^{\R}_n$:] $H^n \subset\R^{n+1}$ is a complete hypersurface in flat Euclidean space $(\R^{n+1},g_{eucl})$ with $0\in H$,
  which is an oriented minimal boundary of some open set in $\R^{n+1}$.
  \item[${\cal{H}}_n$:] ${\cal{H}}_n:= {\cal{H}}^c_n \cup {\cal{H}}^{\R}_n$ and ${\cal{H}} :=\bigcup_{n \ge 1} {\cal{H}}_n$. We briefly refer to $H \in {\cal{H}}$ as an \textbf{area minimizer}.
  \item[$\mathcal{C}_{n}$:] $\mathcal{C}_{n} \subset {\cal{H}}^{\R}_n$ is the space of area minimizing $n$-cones in $\R^{n+1}$ with tip in $0$.
  \item[$\mathcal{SC}_{n}$:] $\mathcal{SC}_n \subset\mathcal{C}_n$ is the subset of cones which are at least singular in $0$.
  \item[$\mathcal{K}_{n-1}$:] For any area minimizing cone $C \subset \R^{n+1}$ with tip $0$, we get the (non-minimizing) minimal hypersurface $S_C:= \p B_1(0) \cap C \subset S^n \subset  \R^{n+1}$
  and we set ${\cal{K}}_{n-1}:= \{ S_C\,| \, C \in {\mathcal{C}_{n}}\}$. We write ${\cal{K}}= \bigcup_{n \ge 1} {\cal{K}}_{n-1}$ for the space of all such hypersurfaces $S_C$.
\end{description}

\begin{remark}
Any current in ${\cal{H}}_n$ can be locally decomposed into (locally disjoint) oriented minimal boundaries of open sets, cf.\ \cite[Appendix A]{L1}. Consequently, we may assume that $H$ is \emph{locally} an oriented boundary of an open subset of $M$.\qed
\end{remark}

\noindent\textbf{Almost Minimizers.} \,
Area minimizers belong to the larger class of \emph{almost minimizers}. Their defining property  is an asymptotically resemblance to honest area minimizers near their singular set. We refer to \cite{T1}, \cite{T2}, \cite{Bo} and \cite{A} for their basic theory. The following definition and result is taken from \cite[Theorem 1]{T1} (see also \cite[Appendix A]{L1} for more background details).\\

Let $\Omega \subset \R^n$ be open and $E \subset \R^n$. Then $\p E$ is called an \emph{almost minimizing boundary} in $\Omega$ if for some $K > 0$, $\alpha \in (0,1)$, the following almost optimal \emph{isoperimetric inequality}
\begin{equation}\label{iso}
\int_{B_\rho(x)}|D \chi_{E}| - \inf \left\{\int_{B_\rho(x)}|D \chi_{F}|\, \Big| \, F
\Delta E \subset \! \subset B_\rho(x) \right\} \le K \cdot \rho^{n-1+2 \cdot \alpha}
\end{equation}
holds for all $x \in \Omega$, $\rho\in(0,dist(x,\p \Omega))$. Here,
\begin{itemize}
  \item $F \Delta E := F \setminus E \cup  E \setminus F$ and $\chi_A$ is the characteristic function of the set $A \subset \R^n$.
  \item $\int_\Omega|D \chi_A| := \sup \{\int_\Omega \chi_A \cdot \mbox{div} g  \, d\mu \;| \; g \in C_0^1(\Omega,\R^n), |g|_{C^0} \le 1 \}$ is called the \emph{perimeter} of $A$ which one may interpret as the area of $\p A$ in $\Omega$.
\end{itemize}
Standard regularity theory says that an almost minimizing boundary is a $C^{1, \alpha}$-hypersurface except for some singular set of Hausdorff-dimension $\le n-8$. Further, any sequence of \emph{blow-ups}, that is, a sequence of infinite rescalings of the almost minimizer, subconverges to a minimal boundary. More generally, we call \emph{almost minimizer} any connected integer multiplicity rectifiable current of dimension $n \ge 2$ which admits a local decomposition into (locally disjoint) almost minimizing boundaries of open sets.

\begin{remark}\label{ram} \,
1.\ In the more general case of a smooth Riemannian manifold $N^m$ as ambient space we call a hypersurface $H^{n-1} \subset N^n$ an almost minimizer if the following condition holds. For each point $p \in H$ there is a ball $B \subset N^n$ centered in $p$, as well as a diffeomorphism $\phi:B\to B_1(0)\subset \R^n$ with  $\phi(p)=0$ such that condition \eqref{iso} holds for $\phi(H \cap B)$ near $0$.\\

2.\ Since the left hand side of \eqref{iso} vanishes, area minimizing hypersurfaces are in particular almost minimizers. On the other hand, the hypersurfaces in ${\cal{K}}$ are almost minimizers which are neither area minimizing nor stable. Condition \eqref{iso} holds for any compact \emph{smooth} hypersurface.\qeda
\end{remark}

We shall consider the following classes of almost minimizers:
\begin{description}
  \item[${\cal{G}}^c_n$:] $H^n \subset M^{n+1}$ is a compact embedded \emph{almost} minimizer with $\p H \v$. We set ${\cal{G}}^c :=\bigcup_{n \ge 1} {\cal{G}}^c_n.$
  \item [${\cal{G}}_n$:] ${\cal{G}}_n := {\cal{G}}^c_n \cup {\cal{H}}^{\R}_n$ and ${\cal{G}} :=\bigcup_{n \ge 1} {\cal{G}}_n$.
\end{description}
${\cal{G}}$ is the main class considered in this paper. Note that even if one is primarily interested in ${\cal{G}}^c_n$ (or ${\cal{H}}^c_n$ ) it is important to add ${\cal{H}}^\R_n$ for this makes the space ${\cal{G}}_n$ closed under blow-ups around singular points. Then we can use compactness results, for the space ${\cal{H}}^\R_n$, to study $H \in {\cal{G}}^c_n$ near $\Sigma_H$. Finally, the theory also extends to the case where $\p H \n$, see Remark~\ref{mgr}.\\

\noindent\textbf{One-Point Compactifications.} \,
As a last piece of notation we introduce the one-point compactification of a hypersurface $H \in {\cal{H}}^{\R}_n$ which we denote by $\widehat{H}$. For the singular set $\Sigma_H$ of $H \in {\cal{H}}^{\R}_n$ we \emph{always} add the point at infinity $\infty_H$ to $\Sigma$ and set $\widehat\Sigma:=\Sigma \cup \infty_H$. This \emph{also} applies in the case where $H$ is non-compact, but $\Sigma$ compact. For $H \in {\cal{G}}^{c}_n$, we set $\widehat H=H$ and $\widehat\Sigma=\Sigma$.

%
%
\subsubsection{\si-Structures}\label{sist}
In this section we review \emph{\si-structures} on hypersurfaces in $\cal{G}$ from \cite{L1}. These are the basic ingredients for our asymptotic analysis near the singular set. \\

\noindent\textbf{\si-Transforms.} \,
An \emph{\si-transform} $\bp_H$ results from merging the metric $g_H$ on $H \setminus \Sigma_H$ and the second fundamental form $A=A_H$ into one scalar function $\bp_H$  on $H \setminus \Sigma_H$ so that $\bp_H$ satisfies the following axioms. The existence of \si-transforms is proved in \cite[Theorem 1.5 and A.III]{L1}.

\begin{definition}[\si-transforms]\label{def1}
We call an assignment $\bp$ that associates with any $H \in {\cal{G}}$ a function $\bp_H:H \setminus \Sigma_H\to\R$ an \textbf{\si-transform} provided it satisfies the following axioms:
\begin{description}
  \item[(S1)] \emph{\textbf{Trivial Gauge}} \,
  If $H \subset M$ is totally geodesic, then $\bp_H \equiv 0$.
  \item[(S2)] \emph{\textbf{\si-Properties}} \,
  If $H$ is not totally geodesic, then the level sets $\M_c:= \bp_H^{-1}(c)$, for $c>0$, we call the $|A|$-\textbf{skins}, surround the level sets of $|A|$:
  \[
  \bp_H>0, \bp_H \ge |A_H| \mm{ and } \bp_H(x) \ra \infty, \mm{ for } x \ra p \in \Sigma_H.
  \]
Like $|A_H|$, $\bp_H$ anticommutes with scalings, i.e., $\bp_{\lambda \cdot H} \equiv \lambda^{-1} \cdot  \bp_{H}$ for any $\lambda >0$.
  \item[(S3)] \emph{\textbf{Lipschitz regularity}} \,
  If $H$ is not totally geodesic, and thus $\bp_H>0$, we define the \textbf{\si-distance}  $\delta_{\bp_H}:=1/\bp_H$.
  This function is
  $L_{\bp}$-Lipschitz regular for some constant $L_{\bp}=L(\bp,n)>0$, i.e.,
  \[
  |\delta_{\bp_H}(p)- \delta_{\bp_H}(q)|   \le L_{\bp} \cdot d_H(p,q) \mm{ for any } p,q \in  H \setminus \Sigma \mm{ and any } H \in {\cal{G}}_n.
  \]
If $H$ is totally geodesic, and thus $\bp_H \equiv 0$, we set $\delta_{\bp_H}\equiv\infty$ and  $|\delta_{\bp_H}(p)- \delta_{\bp_H}(q)|\equiv0$.
  \item[(S4)]  \emph{\textbf{Naturality}} \,
  If $H_i \in {\cal{H}}_n$, $i \ge 1$, is a sequence converging* to the limit space $H_\infty \in {\cal{H}}_n$, then $\bp_{H_i}\overset{C^\alpha}  \longrightarrow {\bp_{H_\infty}}$ for any $\alpha \in (0,1)$. For general $H \in {\cal{G}}_n$, this holds for blow-ups:  $\bp_{\tau_i \cdot H} \overset{C^\alpha}  \longrightarrow {\bp_{H_\infty}}$, for any sequence $\tau_i \ra \infty$ so that  $\tau_i \cdot H \ra H_\infty \in {\cal{H}}^\R_n$.
\end{description}
\end{definition}
*For the precise notions of convergence is explained in \cite[Ch.1.3 and A.III]{L1}. We omit the index $H$ in $\bp_H$ and $\delta_{\bp_H}$ if there is no risk of confusion.

\begin{remark}\label{tg} In this theory the totally geodesic hypersurfaces in  ${\cal{G}}$ play the role of the trivial case. They are always smooth submanifolds, cf. \cite[Corollary A.6]{L1} and in the non-compact case they are just Euclidean hyperplanes. In this case, many results in this paper either obvious or they degenerate to conventions.\qed
\end{remark}

The Lipschitz regular \si-distance $\delta_{\bp}$ admits a Whitney type $C^\infty$-smoothing $\delta_{\bp^*}$ which satisfies (S1)-(S3) and is quasi-natural in the sense that $c_1 \cdot \delta_{\bp}(x) \le \delta_{\bp^*}(x)  \le c_2 \cdot \delta_{\bp}(x)$, for some constants $c_1$, $c_2>0$, cf.\ \cite[Proposition B.3]{L1}.\\

Throughout this paper we choose one fixed \si-transform $\bp$. The precise choice is immaterial for the sequel as the results will not depend on the concrete \si-transform.\\

\noindent\textbf{\si-Uniformity.} \,
As a key application of \si-transforms we can formulate and prove that $\Sigma$ can be approached in a quantitatively non-tangential way from $H \setminus \Sigma$. This can be regarded as a global boundary regularity for $H\setminus\Sigma$.\\

We recall some definitions. A curve will be a continuous map $\gamma:[a,b] \ra X$, $a<b$, into a metric space $(X,d)$. Its \textbf{length} is defined by {\small $l(\gamma):= \sup\{\sum_{i=0,..,N}d(\gamma(t_{i-1}),\gamma(t_i))\,\Big|\, \mm{all partitions }a=t_0 \le t_1 \le ... \le t_N = b\}$}. A curve $\gamma$ is called \textbf{rectifiable} if $l(\gamma) < \infty$. The metric space $X$ is \textbf{rectifiably connected} if each pair of points in $X$ can be joined by a rectifiable curve.

\begin{definition}[Uniform spaces]\label{ud} \, Let $(X,d)$ be a non-complete, locally compact, locally complete and rectifiably connected metric space. We denote its metric completion by $\overline{X}$ and define its boundary by $\p X:= \overline{X} \setminus X$. Then $(X,d)$ is called a \textbf{c-uniform space}, or a \textbf{uniform space} for short, if there is some $c\ge 1$ such that any two points $p,q \in X$ can be joined by a \textbf{c-uniform curve}. This is a rectifiable curve $\gamma_{p,q}: [a,b] \ra X$ going from $p$ to $q$ such that:
\begin{itemize}
  \item \emph{\textbf{Quasi-geodesic:}} \, $l(\gamma_{p,q})\le c\cdot d(p,q)$.
  \item \emph{\textbf{Twisted double cones:}} \, For any $z \in \gamma_{p,q}$ let $l_{min}(\gamma_{p,q}(z))$ be the minimum of the lengths of the two subcurves of $\gamma_{p,q}$ from $p$ to $z$ and from $z$ to $q$. Then
  \[
  l_{min}(\gamma_{p,q}(z)) \le c \cdot dist(z,\p X).
  \]
\end{itemize}
\end{definition}

For instance, bounded domains in $\R^n$ with Lipschitz regular boundary are uniform, and so are fractal like spaces such as snowflakes or the complement of the Sierpinski gasket~\cite{Ai1}, \cite{ALM}. On the other hand,  domains such as the cube $C= (-1,1)^n  \subset \R^n$ is no longer uniform after deleting the inscribed ball $B=B_1(0)  \subset \R^n$, $n \ge 3$ - the uniformity of $C \setminus B$ is violated near $\p C \cap \p B$.\\

Now for any $H \in {\cal{G}}_n$ the regular part $H \setminus \Sigma$ is uniform in an actually sharper sense which also takes the curvature of $H \setminus \Sigma$ into account. This is the \emph{\si-uniformity} of $H \setminus \Sigma$, cf.\ \cite[Theorem 1.8 and Appendix A.III]{L1}.

\begin{theorem}[\si-Uniformity of $\mathbf{H \setminus \Sigma}$]\label{thm2} \,
For any $H \in {\cal{G}}_n$ with (possibly empty) singular set $\Sigma=\Sigma_H$ holds the following.
\begin{enumerate}
  \item$H \setminus \Sigma$ and $H$ are \textbf{rectifiably connected}. In particular, any $H \in {\cal{G}}^c_n$ has a finite intrinsic diameter: $diam_{g_H}H <\infty.$
  \item There exists $c>0$ such that $H \setminus \Sigma$ is a \textbf{c-\si-uniform space}, or \textbf{\si-uniform space} for short. This means that any pair $p,q \in H \setminus \Sigma$ can be joined by a \textbf{c-\si-uniform curve} in $H \setminus \Sigma$, i.e., a rectifiable curve $\gamma_{p,q}: [a,b] \ra H \setminus \Sigma$ with $\gamma_{p,q}(a)=p$, $\gamma_{p,q}(b)=q$ and so that:
  \begin{itemize}
    \item \emph{\textbf{Quasi-geodesic:}} \, $l(\gamma)  \le c \cdot  d(p,q).$
    \item \emph{\textbf{Twisted double \si-cones:}} \, $l_{min}(\gamma_{p,q}(z)) \le c \cdot \delta_{\bp}(z)$ for any $z \in \gamma_{p,q}$.
  \end{itemize}
  \item For $H\in{\cal{H}}^{\R}_n$ the \si-uniformity parameter $c$ depends only on $n$. Moreover for any compact family in $H \in{\cal{H}}^c_n$ we can choose a common \si-uniformity parameter.
  \item For $H\in{\cal{G}}_n$ the rescalings $k \cdot H$, for integers $k\ge 1$, subconverge to a tangent cone $C\in{\cal{H}}^{\R}_n$. Then there exists a common $c_H>0$ such that the $k \cdot H$ and $C$ are $c_H$-\si-uniform.
\end{enumerate}
\end{theorem}

\begin{remark}
The Lipschitz condition (S3) implies $\delta_{\bp}(x) \le L  \cdot dist(x,\Sigma)$ for any $x \in H \setminus \Sigma$. Thus, if $\Sigma \n$, \si-uniformity of $H \setminus \Sigma$ implies uniformity -- a result which would be hard to derive directly. Note, however, that \si-uniformity still makes sense in the regular case where $\Sigma \v$.\qed
\end{remark}

%
%
\subsubsection{Main Results of the Present Paper}\label{pre}
The analysis of elliptic operators on a manifold near its boundary is controlled by elliptic boundary regularity results such as, most basically, boundary Harnack inequalities. The control and the inequalities depend both on properties of the manifold and of the operator towards but, less obviously, also away from the boundary. For the Laplace operator $\Delta$ on $\R^n$, boundary Harnack inequalities hold on any uniform domain $D\subset \R^n$. Concretely, we say that the pair $(D,\Delta)$ satisfies the \emph{boundary Harnack principle} if the subsequent property is satisfied:\\

\noindent\textbf{Boundary Harnack Principle (BHP):} \emph{There exist constants $A$, $C>1$ depending only on $D \subset \R^n$ such that for any point $p\in \p D$ and sufficiently small $R>0$ the following is true. For any two harmonic functions $u$, $v>0$ on $B_{A \cdot R}(p) \cap D$ which vanish along $B_{A \cdot R}(p) \cap \p D$},
\begin{equation}\label{gf2}
u(x)/v(x) \le C \cdot  u(y)/v(y) \mm{ for all } x,\, y \in B_R(p) \cap D.\\
\end{equation}

By a result due to Aikawa~\cite{Ai1}, \cite{Ai2}, \cite{Ai3}, the \emph{uniformity} of $D$ is actually the minimal inner regularity condition needed to ensure the validity of the BHP.

\begin{remark}
In the literature, for instance \cite{Ai1}, this formulation of the BHP is also called the \emph{local} BHP. It is this condition which garantuees that the Martin boundary is homeomorphic to the topological boundary, see below. A weaker version is the \emph{global} BHP. It merely asserts that there is a Harnack constant $C$ which depends also on $p$ and $R$.\qed
\end{remark}

\noindent\textbf{BHP and Martin Theory on $\mathbf{H \setminus \Sigma}$.}
Next we pass from uniform domains in Euclidean space to area minimizers. However, the uniformity condition discussed above does not take the curved and degenerating geometry of $(H \setminus \Sigma, g_H)$ into account. It is precisely the stronger \si-uniformity which handles these additional geometric difficulties. To describe the elliptic problems on $(H \setminus \Sigma, g_H)$ that we can address this way we use special charts for $H \setminus \Sigma$, namely \emph{\si-adapted charts}. These are certain bi-Lipschitz charts $\psi_p:B_R(p)\to\R$ centered in $p\in H\setminus\Sigma$ where the radius $R$ of the ball depends in particular on $\bp_H(p)$, cf.\ \cite[Proposition B.1]{L1} and Chapter \ref{agrm} below.

\begin{definition}\label{sao}
Let $H \in \cal{G}$. A second order elliptic operator $L$ on $H\setminus \Sigma$ is called \textbf{\si-adapted} if the following two conditions hold.
\begin{itemize}
  \item \emph{\textbf{$\mathbf{\bp}$-adaptedness:}}
  There exists a constant $k=k_L\ge 1$ such that $L$ satisfies the following \si-weighted uniformity conditions. With respect to the charts $\psi_p$, $L$ can be locally written
  {\small \[
  -L(u) = \sum_{i,j}  a_{ij} \cdot \frac{\p^2 u}{\p x_i \p x_j} + \sum_i b_i \cdot \frac{\p u}{\p x_i} + c \cdot u
  \]}
  for $\beta$-H\"{o}lder continuous coefficients $a_{ij}$, $\beta \in (0,1]$, and measurable functions $b_i$, $c$, such that for all $p\in H\setminus\Sigma$, $\xi\in\R^n$, we have
  \begin{enumerate}
    \item $k^{-1} \cdot\sum_i \xi_i^2 \le \sum_{i,j} a_{ij}(p) \cdot \xi_i \xi_j \le k \cdot \sum_i \xi_i^2$,
    \item $\delta^{\beta}_{\bp}(p) \cdot  |a_{ij}|_{C^\beta(B_{\theta(p)}(p))} \le k$,  $\delta_{\bp}(p) \cdot |b_i|_{L^\infty} \le k$ and $\delta^2_{\bp}(p) \cdot |c|_{L^\infty} \le k$.
  \end{enumerate}
  \item \emph{\textbf{$\mathbf{\bp}$-weak coercivity:}
  There exists a positive $C^2$-supersolution $u$ of $L \, f = 0$ with}
  \[
  L \, u \ge \ve \cdot \bp^2 \cdot u\mm{ for  some } \ve  >0.
  \]
\end{itemize}
\end{definition}
We denote the largest such $\ve$ by $\ve_L$. In the case of symmetric operators we will see, cf.\ Theorem~7 below, that $\ve_L$ is the \emph{(generalized) principal eigenvalue} of $\delta_{\bp}^2  \cdot L$.

\begin{remark} 1.\ Provided that $L$ has sufficiently regular coefficients, it is enough to consider weak (super)solutions to infer the existence of regular (super)solutions.\\

2.\ The $\bp$-weak coercivity guarantees, in particular, the existence of a \emph{Green's function} $G:H\setminus\Sigma\times H\setminus\Sigma\to(0,\infty]$. That is, $G(\cdot,y)$ solves $L\,f=0$ on $H\setminus(\Sigma\cup\{y\})$ and $L\,G(\cdot,y)= \delta_y$, the Dirac measure in $y$, cf.\ \cite[Paragraph 1.2]{An1}. In the sequel, $G(x,y)$ will denote the \emph{minimal} Green's function of an \si-adapted operator $L$, see the beginning of Section~\ref{nmt}.\\

3.\ In \cite[Theorem 2]{L2} we will see that, independently of the chosen \si-transform $\bp$, many classical elliptic operators are actually \si-adapted.\qeda
\end{remark}

On $H \setminus \Sigma$ a solution $u >0$ of the equation $L \, f=0$ will usually diverge to infinity when we approach $\Sigma$, and so will a Green's function of $L$. The  generalization of the vanishing boundary data condition for the Laplacian on a Euclidean domain is a \emph{minimal growth} condition: A solution $u >0$ of $L \, f=0$ is said to be \textbf{\emph{L}-vanishing} in $p \in \widehat{\Sigma}$, the one-point compactifications of $\Sigma$ cf.~\ref{notation}, if there is a supersolution $w >0$,   such that $u/w(x) \ra 0$, for $x \ra p$, $x \in H \setminus \Sigma$.  For $H \in {\cal{G}}$ we show that \emph{L}-vanishing solutions satisfy a BHP, where $\widehat{\Sigma}$ plays the role of the  boundary.\\

To formulate our boundary Harnack inequalities on $H\setminus \Sigma$, for an $H \in {\cal{G}}$, we replace the systems of metric balls in the Euclidean BHP of (\ref{gf2}) with  some particular type of neighborhood  basis $\mathbf{N}^\delta_i(z) \subset H$, $i \in \Z^{\ge 0}$, of points $z \in \widehat{\Sigma}$ with  $\mathbf{N}^\delta_{i+1}(z)  \subset \mathbf{N}^\delta_i(z)$ and $\bigcap_i \mathbf{N}^\delta_i(z)=\{z\}$. We postpone the details to Ch.\ref{bhai}. For now we note that a typical choice for the subsets $\mathcal{N}^\delta_i(z) :=\mathbf{N}^\delta_i(z) \cap H \setminus \Sigma$ are halfspaces in the hyperbolic \si-geometry on $H\setminus \Sigma$ introduced in~\cite[Theorem 1.11]{L1}.  This can be motivated from the Poincar\'{e} disc model for the hyperbolic plane where Euclidean discs (around the boundary points of the unit disc $D \subset \R^{2}$) become hyperbolic halfspaces.\\

\noindent\textbf{Theorem 1} \textbf{(BHP on ${\cal{G}}$, see Theorem~\ref{mbhsq})} \,
\emph{Let $H \in {\cal{G}}$ and $L$ be an \si-adapted operator on $H\setminus \Sigma$. Then there exists a constant $C(H,L) >1$ such that for any $z \in \widehat{\Sigma}$ and any  two solutions
 $u$, $v >0$ of $L\, f= 0$ on $H \setminus \Sigma$ which are both L-vanishing along $\mathbf{N}^\delta_i(z)\cap \widehat{\Sigma}$, we have}
\begin{equation}\label{fhepq1}
u(x)/v(x) \le C \cdot  u(y)/v(y), \mm{ \emph{for all} }x,\, y \in \mathcal{N}^\delta_{i+1}(z) \mm{ and for any $i \in \Z^{\ge 1}$}.
\end{equation}

The constant $C$ only depends on rather coarse data we can extract from $H \subset M$ and $L$. The case of Euclidean area minimizers is particularly neat:\\

\noindent\textbf{Theorem 2} \textbf{(Stable BHP on ${\cal{H}}^{\R}_n$, see Theorem~\ref{mbhsqq})} \,
\emph{For any $\kappa$, $\eta>0$, there exists a constant $C(\kappa, \eta, n)>1$ so that for any $H \in {\cal{H}}^{\R}_n$ and any  \si-adapted operator $L$ on $H \setminus \Sigma$ the following stable form of a BHP holds:}\\

\emph{For $k_L \le \kappa$, $\ve_L \ge \eta$  and $z \in \widehat{\Sigma}$ and any two solutions $u$, $v >0$ of $L \, f= 0$ on $H \setminus \Sigma$ which are both L-vanishing along $\mathbf{N}_i(z)\cap \widehat{\Sigma}$, we have}
\begin{equation}\label{fhepq1a}
u(x)/v(x) \le C \cdot  u(y)/v(y), \mm{ \emph{for all} } x, y \in \mathcal{N}^\delta_{i+1}(z)   \mm{ and for any $i \in \Z^{\ge 1}$}.
\end{equation}

Boundary Harnack inequalities are the fundamental not only in the context of local boundary regularity results. They can be used more globally to characterize the extremal elements amongst the solutions of $L \, f=0$: the minimal solutions. We call $u >0$ \textbf{minimal} if for any other solution $v >0$ with $v \le u$, we have $v \equiv c \cdot u$ for some $c >0$.  The space of minimal solutions (normalized to $1$ in some basepoint) is fundamental in potential theory. It is called the \textbf{(minimal) Martin boundary}  $\p^0_M (H \setminus \Sigma,L)$.   A priori, $\p^0_M(X,L)$ may heavily depend on both the operator $L$ and the underlying space $X$. Even for the Laplace operator on rather symmetric spaces this boundary is generally hard to describe cf. \cite{GLT}.  In view of this the following result is surprising. (We recall the mentioned notions from Martin theory  in Ch.\ref{nmt} below.) \\

\noindent\textbf{Theorem 3} \textbf{(Martin Theory on $\mathbf{H \setminus \Sigma}$, see Theorem~\ref{mbhs})} \,
\emph{Let $H \in {\cal{G}}$ be a non-totally geodesic hypersurface and $L$ some \si-adapted operator on $H \setminus \Sigma$. Then
\begin{itemize}
  \item the identity map on $H \setminus \Sigma$ extends to a homeomorphism  between $\widehat{H}$ and the \emph{\textbf{Martin compactification}} $\overline{H \setminus \Sigma}_M$.
  \item all Martin boundary points are minimal:
  $\p^0_M (H \setminus \Sigma,L) \equiv \p_M(H \setminus \Sigma,L)$.
\end{itemize}
In particular, $\widehat{\Sigma}$ and the minimal Martin boundary $\p^0_M (H \setminus \Sigma,L)$ are homeomorphic. \\}

Using the Choquet integral representation in \cite[Chapter 6]{C} we obtain from Theorem 3 the following general version of the Martin representation theorem. A function $u >0$ on $H \setminus \Sigma$ solves $L \, f= 0$ if and only if there is a unique finite Radon measure $\mu$ on $\widehat{\Sigma}$ such that
\begin{equation}\label{muu}
u(x) = u_{\mu}(x) =\int_{\widehat{\Sigma}} k(x;y) \, d \mu(y).
\end{equation}
In this representation formula, $k(x;y)$ denotes the \emph{Martin kernel} of $L$ on $H \setminus \Sigma$. It is, up to multiples, the unique positive solution of $L \, f = 0$ on $H \setminus \Sigma$ which  $L$-vanishes in all points of $\widehat{\Sigma}$ except for $y$. Moreover, the functions $k(\cdot;y)$, $y\in\p_M(H \setminus \Sigma,L)$, are just the minimal solutions of $L$.\\

\noindent\textbf{Extension Results to $\mathbf{\Sigma}$.} \,
The following three extension theorems can be regarded as counterparts of classical results for the Laplacian on suitably regular Euclidean domains (e.g., uniform domains), cf.\ for instance \cite{AG}, \cite{Ai2} or \cite{JK}. We start with a Fatou type theorem for \si-adapted operators on $H \setminus \Sigma$. For this, we quantify the non-tangentiality of $\Sigma$ thought of as a boundary by means of non-tangential \si-pencils
\[
\P(z,\omega):= \{x \in H \setminus \Sigma \,|\, \delta_{\bp} (x) >\omega \cdot d_H(x,z)\}
\]
pointing to $z \in \Sigma$, where $\omega>0$. One may think of the angle $\arctan(\omega^{-1})$ as some kind of aperture of $\P(z,\omega)$ relative to $z$.\\

\noindent\textbf{Theorem 4} \textbf{(Relative Fatou Theorem on  $\mathbf{H \setminus \Sigma}$, see Theorem~\ref{ftt})} \,
\emph{Let $H \in {\cal{G}}$ and $L$ be an \si-adapted operator on $H \setminus \Sigma$. Further, let $\mu$ and $\nu$ be two finite Radon measures on $\Sigma$ associated with solutions $u_{\mu}$ and $u_{\nu}$ of $L \, f=0$, cf.\ \eqref{muu}. Then for $\nu$-almost any $z \in \Sigma$ and any fixed $\omega >0$, we have
\[
u_{\mu}/u_{\nu}(x) \ra d \mu/d \nu(z)\mm{ as } x \ra z, \mm{ with } x \in \P(z,\omega).
\]
(Here, $d \mu/d \nu$ denotes the Radon-Nikodym derivative of $\mu$ with respect to $\nu$.)}\\

In the case of an open subset $A \subset \Sigma$ with $u_{\mu}(A)=u_{\nu}(A)=0$ the Fatou theorem does not give any information on the behavior of
$u_{\mu}/u_{\nu}(x)$ as $x \ra z \in A$. However, we can still invoke boundary Harnack inequalities since in this case, both $u_{\mu}$ and $u_{\nu}$ \emph{L}-vanish along $A$. We obtain the following complementary result which even applies tangentially.\\

\noindent\textbf{Theorem 5} \textbf{(Continuous Extensions to $\mathbf{\Sigma}$, see Theorem~\ref{cee})} \,
\emph{Let $H \in {\cal{G}}$ and $L$ be an \si-adapted operator on $H \setminus \Sigma$. Then for any two solutions  $u$, $v >0$ of $L\, f = 0$ on $H \setminus \Sigma$ which are both  $L$-vanishing along some common open set $A \subset\widehat{\Sigma}$, the quotient $u/v$ on $H \setminus \Sigma$ admits a continuous extension to $(H \setminus \Sigma) \cup A$.}\\

Unlike their quotients, the individual solutions $u$ and $v$ usually diverge towards $\Sigma$. To attack boundary value problems we need, in addition to \si-adaptedness, some a priori control over certain solutions as for instance in the following result.\\

\noindent\textbf{Theorem 6} \textbf{(Dirichlet Problem for \si-Adapted Operators, see Theorem~\ref{diri})} \,
\emph{Let $H \in {\cal{G}}$ and $L$ be an \si-adapted operator on $H \setminus \Sigma$ such that
\begin{itemize}
  \item constant functions solve $L \, f=0$.
  \item for a given $p \in H \setminus \Sigma$, the Green's function $G(x,p)\ra 0$ as $x \ra \widehat{\Sigma}$.
\end{itemize}
Then, for any continuous function $f$ on $\widehat{\Sigma}$, there exists a uniquely determined continuous function $F$ on $H$ such that}
\[
L \, F =0 \mm{ with } F|_{\widehat{\Sigma}} \equiv f.
\]
\

\noindent\textbf{Symmetric Operators.} \,
A frequently considered type of elliptic problems is that of eigenvalues. For this we shall focus on symmetric operators where it is possible and useful to bring the weak coercivity condition into a variational form:\\

\noindent\textbf{Theorem 7} \textbf{(\si-Hardy Inequality, see Theorem~\ref{hi})} \,
\emph{Let $H \in{\cal{G}}$ and $L$ be a \emph{symmetric}, $\bp$-adapted operator on $H\setminus \Sigma$. Then the weak coercivity condition  is equivalent to the existence of a positive constant $\tau = \tau(L,\bp,H)>0$  such that the Hardy type inequality
\begin{equation}\label{hadi0}
\int_H  f  \cdot  L f  \,  dV \, \ge \, \tau \cdot \int_H \bp^2\cdot f^2 dV
\end{equation}
holds for any smooth $f$ which is compactly supported in $H\setminus \Sigma$.}\\

For a given $\bp$-adapted symmetric operator $L$ there exists a largest $\lambda^{\bp}_{L,H} \in [-\infty,+\infty)$ such that the Hardy inequality \eqref{hadi0} is satisfied. This constant can be viewed as an eigenvalue for the operator $\delta_{\bp}^2 \cdot L$, the so-called \emph{(generalized) principal eigenvalue} of $\delta_{\bp}^2 \cdot L$, cf.\ \cite[Chapter 4]{P} for the associated spectral geometry on unbounded domains. We notice that $L$ is \si-adapted if $\lambda^{\bp}_{L,H} >0$. We use this simple observation to extend Definition \ref{sao}:

\begin{definition}
An $\bp$-adapted symmetric operator $L$ on $H \setminus \Sigma$ is called \textbf{shifted \si-adapted} if the principal eigenvalue of $\delta_{\bp}^2 \cdot L$ is finite, i.e., $\lambda^{\bp}_{L,H} >-\infty$.
\end{definition}

 The r\^ole of the principal eigenvalue is explained by the following trichotomy, which resembles that for the spectral theory of operators on general Euclidean domains, cf.\ \cite[Chapter 4]{P}.\\

\noindent\textbf{Theorem 8} \textbf{(Criticality, see Theorem~\ref{scal2})} \,
\emph{For any singular $H \in \cal{G}$ and any shifted \si-adapted operator $L$ on $H \setminus \Sigma$ with H\"older continuous coefficients we set
\[
L_\lambda:= L - \lambda \cdot \bp^2 \cdot Id,\,\,\lambda \in\R.
\]
Then we have the following trichotomy.
\begin{itemize}
    \item \emph{\textbf{Subcritical case}} when $\lambda < \lambda^{\bp}_{L,H}$: The operator $L_\lambda$ is \si-adapted, and the minimal solutions of $L_\lambda \, v=0$ $L$-vanish in all but the one point in $\widehat{\Sigma}$ which represents this solution as a Martin boundary point.
    \item \emph{\textbf{Critical case}} when $\lambda = \lambda^{\bp}_{L,H}$: There is, up to multiples, a
    \[
    \mm{unique positive solution, the so-called \textbf{ground state} } \phi  \mm{ of L} _{\lambda^{\bp}_{L,H}} \, f= 0.
    \]
    The ground state $L$-vanishes along $\widehat{\Sigma}$ and can be described as the limit of first Dirichlet eigenfunctions for the operator $\delta_{\bp}^2 \cdot L$ on a sequence of smoothly bounded domains $\overline{D}_m \subset D_{m+1} \subset H \setminus \Sigma$, $m \ge 0$, with $\bigcup_m D_m = H \setminus \Sigma$.
    \item \emph{\textbf{Supercritical case}} when $\lambda > \lambda^{\bp}_{L,H}$: $L_\lambda \, f = 0$ has \emph{no} positive solution.
\end{itemize}
}

\begin{remark}
Each of these theorems admits an extension to the case of (almost) minimizers with non-empty boundary $\p H \n$ where the operators may also degenerate towards $\p H$, see Remark \ref{mgr}.\qed
\end{remark}

%
%
%
\setcounter{section}{2}
\renewcommand{\thesubsection}{\thesection}
\subsection{Ancona's Hyperbolic Boundary Harnack Principle}
%
Ancona \cite{An1, An2} developed a potential theory applicable to many elliptic operators on Gromov hyperbolic manifolds of bounded geometry.
 In this chapter we briefly review some essential concepts. In \cite{KL} we have given a detailed account on this theory and we oftentimes use it
 as a reference for some concrete results not explicitly stated in other sources.
%
\subsubsection{Gromov Hyperbolic Manifolds and $\Phi$-Chains}\label{app}

We first recall the notion of Gromov hyperbolicity.  It has no local impact but strong consequences for the geometry near infinity. We mention \cite{BH} and \cite{BHK} as general references and \cite{KL} for some complementary discussions.

\begin{definition}[Gromov Hyperbolicity and Gromov Boundary]\label{del}
We call a geodesic metric space $X$ \textbf{Gromov hyperbolic}, or more precisely $\mathbf{\delta}$\textbf{-hyperbolic,} if its geodesic triangles are $\mathbf{\delta}$\textbf{-thin} for some $\delta=\delta_X >0$. That is, each point on the edge of any geodesic triangle is within $\delta$-distance of one of the other two edges.\\

Two rays in $X$ are \emph{equivalent} if they have finite Hausdorff distance. The set $\p_G X$ of equivalence classes $[\gamma]$ of geodesic rays from a fixed base point $p \in X$ is called the \textbf{Gromov boundary} of $X$. This definition of $\p_G X$ is independent of $p$, cf.\ \cite[Part III.H]{BH}.
\end{definition}

We define a (metrizable) topology on $\overline{X}_G = X \cup \p_G X$ as follows. A \emph{generalized geodesic ray} $\gamma: I \ra X$ is an isometric embedding of the interval $I \subset \R$ into $X$, where either $I = [0,\infty)$, or $I = [0,R]$ for some $R \in (0,\infty)$. In the former case of an infinite interval we call $\gamma$ a \emph{proper geodesic ray}. We fix a base point $p \in X$ and use the hyperbolicity to canonically identify any $x\in X$ with the generalized ray $\gamma_x$ with endpoint $\gamma(R) = x$. If $R<\infty$ we extend the domain to $[0,\infty]$ by setting $\gamma(t) := \gamma(R)$, $t \in [R,\infty]$.

\begin{definition}[Gromov Compactification]
The topology on $\overline{X}_G$ is induced from the following notion of convergence: A sequence $x_n \in \overline{X}$ converges to $x \in \overline{X}$ if there exist generalized rays $\gamma_n$ with $\gamma_n(0) = p$ and $\gamma_n(\infty) = x_n$ subconverging on compact sets to a generalized ray $\gamma$ with $\gamma(0) = p$ and $\gamma(\infty)=x$. It is easy to show, cf.\ \cite{BH}, that
\begin{itemize}
  \item the canonical map $X \hookrightarrow \overline{X}_G$ is a homeomorphism onto its image,
  \item $\overline{X}_G$ is a compact metrizable space,
  \item $\p_G X$ is closed and thus a compact subset of $\overline{X}_G$.
\end{itemize}
The space $\overline{X}_G$ is called the \textbf{Gromov compactification} of $X$.
\end{definition}

To formulate the hyperbolic BHP we use hyperbolic counterparts of the concentric balls used in the classical Euclidean BHP of (\ref{gf2}). They are defined using $\Phi$-chains a concept introduced by Ancona~\cite[Definition 5.1, p.\ 93]{An1}. Our slightly modified definition is taken for   ~\cite{KL} to better match related work in \cite[Ch.8]{BHK}. The typical  $\Phi$-chains in Gromov hyperbolic spaces are families of nested halfspaces contracting to one point in the Gromov boundary.

\begin{definition} [$\Phi$-Chains]
\label{def:phi} For a monotonically increasing function $\Phi:[0,\infty)\to(0,\infty)$ with $\Phi(t)\overset{t\to\infty}{\longrightarrow}\infty$, a \textbf{$\Phi$-chain} on a manifold $X$ is a finite or infinite sequence $U_1\supset U_2\supset\cdots\supset U_m$ of open subsets of $X$ together with a sequence of \textbf{track points} $x_1,x_2,\dots,x_m$ so that
\[ x_i \in \p U_i \, , \, \Phi(0)\le d(x_i,x_{i+1})\le 3 \cdot \Phi(0)\,\mm{ and }\, d(x,\p U_{i \pm 1})\ge\Phi(d(x,x_i)),\mm{ for every }x\in \p U_{i}.\]
\end{definition}

The existence of infinite  $\Phi$-chains can be considered as a partial hyperbolicity property of the underlying space. However, classical hyperbolicity, in the sense of constant negative sectional curvature, is \emph{not} yet sufficient to ensure the existence of infinite $\Phi_\delta$-chains. Simple counterexamples are $\Z^2$-coverings of genus $\ge 2$ Riemann surfaces equipped with some hyperbolic metric. They roughly look like the Euclidean plane with $\Z^2$-periodically attached handles. In turn, \emph{Gromov hyperbolic} spaces carry a myriad of infinite $\Phi$-chains. One can build them starting from some geodesic
 $\gamma: (0,  a)  \ra  X$, for $a>0$, i.e., an isometric embedding of $(0,  a)$ in  $X$. Then we define
\begin{equation}\label{pch}
U^\gamma_t:=\{x \in X\,|\, dist\big(x,\gamma([t,a))\big)<dist\big(x,\gamma((0,t]) \big)\}.
\end{equation}

From some explicit computations one gets the following result we cite from \cite[Lemma 8.3 - Proposition 8.10]{BHK}

\begin{lemma}  [Canonical $\Phi_\delta$-Chains on Hyperbolic Spaces]\label{nn}  \label{cpc} For a  $\delta$-hyperbolic space $X$ and a geodesic $\gamma: (0,  a)  \ra  X$, for $a >10^3 \cdot \delta$. We choose  track points $x_i=\gamma(t_i)$, $i=0,...,m+1$, with $x_0=\gamma(0), x_{m+1}=\gamma(a)$, with $d(x_i,x_{i+1})=300\cdot \delta$, for $i <m$, and $d(x_m,x_{m+1})\le300\cdot \delta$. Then there  is some $\Phi_\delta(t)=a_\delta + b_\delta \cdot t$, for $a_\delta, b_\delta  > 0$ depending only on $\delta$, so that
\begin{equation}\label{vu}
\mathcal{N}^\delta_i(\gamma):=U_{t_i} \mm{ form a } \Phi_\delta\mm{-chain with track points }x_i , i=1,...,m.
\end{equation}

We call the $\mathcal{N}^\delta_i(\gamma)$ a \textbf{canonical $\Phi_\delta$-chain} and we also define the open sets $\mathbf{N}^\delta_i$ in the Gromov compactification $X_G$ of $X$ naturally extending $\mathcal{N}^\delta_i(\gamma) \subset X$.
\begin{equation}\label{vvu}
\mathbf{N}^\delta_i :=\mathcal{N}^\delta_i \cup \{z \in \p_G X\,|\, z \mm{ can be represented by a sequence in }\mathbf{N}^\delta_i\} \subset X_G.
\end{equation}
Then, if $\gamma$ is a ray representing some $z \in \p_G X$, the $\mathbf{N}^\delta_i(\gamma)$ are a neighborhood basis  of $z$ in $X_G$.
\end{lemma}

\begin{remark}\label{cgc2}
The reason for the rather large distance $300 \cdot \delta$ between two track points is to find a point, some kind of a hub, between the any two track points so that there are further controlled auxiliary $\Phi_\delta$-chains chosen to link any two given points in  $\p \mathbf{N}^\delta_i(\gamma)$ and $\p \mathbf{N}^\delta_{i+1}(\gamma)$ path through the hub. This is needed in the analytic application of canonical $\Phi_\delta$-chains.\qed
\end{remark}

Besides Gromov hyperbolicity there is a second crucial prerequisite we need for a transparent potential theory. The underlying space must have \emph{bounded geometry}. That is, it is supposed to locally look the same around any point, up to a uniformly controlled deviation. (In the literature the bounded geometry constraint is occasionally ignored albeit it is the more basic condition. The hyperbolicity is used to improve the results we already get for manifolds of bounded geometry.)

\begin{definition}[Bounded Geometry]\label{bog}
We say that an (at least Lipschitz regular) Riemannian manifold $M$ is of $(\varrho,\ell)$-\textbf{bounded geometry} if there exist global constants $\varrho={\varrho_M}>0$ and $\ell=\ell_M \ge 1$ for $M$ such that for each ball $B_{\varrho}(p) \subset M$ there is a smooth $\ell$-bi-Lipschitz chart $\phi_p$ onto an open set $U_p\subset\R^n$ with its Euclidean metric.
\end{definition}

In \ref{smbg} we will see how this basic form of bounded geometry corresponds with other notions for smooth manifolds.

%
\subsubsection{Singularities as Gromov Boundaries}\label{bhai}
To each $H \in {\cal{G}}$ we can assign its \textbf{\si-metric} $d_{\bp}=d_{\bp_H}$ defined by
\begin{equation}\label{sim}
d_{\bp}(x,y) := \inf \Bigl  \{\int_\gamma  \bp \, \, \Big| \, \gamma   \subset  H \setminus \Sigma\mbox{ rectifiable curve joining }  x \mbox{ and } y  \Bigr \}
\end{equation}
for $x$, $y\in H \setminus \Sigma$. The metric $d_{\bp}$ is also well-defined for \emph{smooth} $H$ where $\Sigma\v$. Alternatively, the \si-metric can be written $\bp^2 \cdot g_H$, but this is not a regular Riemannian metric since $\bp$ is merely a locally Lipschitz function. However, there is a Whitney type smoothing process, cf.\ \cite[Proposition B.3]{L1}:

\begin{proposition}[\si-Whitney smoothings]\label{smsk}
For any \si-transform $\bp$ there exists a smoothing $\bp^*$, i.e., a family of smooth functions $\bp_H^*$ defined on $H \setminus \Sigma$ for any $H \in {\cal{G}}$. The smoothing $\bp^*$ still satisfies axioms (S1) - (S3) for \si-transforms while for (S4), we have the inequalities
\begin{equation}\label{smot}
c_1 \cdot \delta_{\bp}(x) \le \delta_{\bp^*}(x)  \le c_2 \cdot \delta_{\bp}(x)\quad\mm{and}\quad |\p^\beta \delta_{\bp^*}  / \p x^\beta |(x) \le c_3(\beta) \cdot \delta_{\bp}^{1-|\beta|}(x)
\end{equation}
for constants $c_i(L_{\bp},H,\beta) >0$, $i=1,\,2,\,3$. Here, $\beta$ is a multi-index for derivatives with respect to normal coordinates around $x \in H \setminus \Sigma$. For $H \in {\cal{H}}^{\R}_n$ we even have uniform constants $c_i(L_{\bp},n,\beta)$.
\end{proposition}

We interpret \eqref{smot} as a weak $\bp^*$-version of the naturality axiom (S4). We therefore refer to \eqref{smot} as \emph{quasi-naturality}. As a counterpart of $d_{\bp}$ we have the
metric $d_{\bp^*}$ which corresponds to the smooth Riemannian manifold $(H \setminus \Sigma, (\bp^*)^2 \cdot g_H)$. For totally geodesic $H \in {\cal{G}}$, both $(H \setminus \Sigma, d_{\bp})$ and $(H \setminus \Sigma, d_{\bp^*})$ are well-defined but they are one-point spaces since $\bp \equiv 0$. Otherwise, $\bp>0$ and, hence,  $(H \setminus \Sigma, d_{\bp})$ and $(H \setminus \Sigma, d_{\bp^*})$ are homeomorphic  to $(H \setminus \Sigma, g_H)$.\\

The \si-metrics have two important properties, we formally state in \ref{thm3} below. They are Gromov hyperbolic and they have bounded geometry. We recall

\begin{remark}\label{smbg}
The Whitney type smoothing process can equally be employed to upgrade the Lipschitz form of bounded geometry, defined in \ref{bog}, to a \emph{bounded geometry of order $k$}, cf.\ \cite[Proposition B.3]{L1} and \cite[Chapter VI.2]{St} for details. This means that the injectivity radius is bounded below by $10 \cdot \varrho$, and the covariant derivatives $\nabla^i R$ of the curvature tensor $R$ are bounded up to order $k$, i.e., $|\nabla^i R| \le \varrho^{-1}, \mm{ for any } i \le k$.\qed
\end{remark}

Now we can formulate the following hyperbolization results for $H \in {\cal{G}}$, cf.\ \cite[Theorem 1.11, Proposition 3.10 and Theorem 1.13]{L1}.

\begin{theorem}[Conformal Hyperbolic Unfoldings]\label{thm3}
For any non-totally geodesic $H \in {\cal{G}}$, the \si-metric $d_{\bp}$ has the following properties:
\begin{itemize}
  \item The metric space $(H \setminus \Sigma, d_{\bp})$ and its \textbf{quasi-isometric Whitney smoothing}, i.e., the smooth Riemannian manifold $(H \setminus \Sigma, d_{\bp^*}) = (H \setminus \Sigma, 1/\delta_{\bp^*}^2 \cdot g_H)$, are \textbf{complete Gromov hyperbolic spaces} with \textbf{bounded geometry}.
  \item $d_{\bp}$ is \textbf{natural}, that is, the assignment $H\mapsto d_{\bp_H}$ commutes with compact convergence of regular domains of the underlying area minimizers.
  \item For any \emph{singular} $H \in {\cal{G}}$ the identity map on $H \setminus \Sigma$ extends to \textbf{homeomorphisms}
  \[
  \widehat{H}\cong\overline{(H \setminus \Sigma,d_{\bp})}_G \cong \overline{(H \setminus \Sigma,d_{\bp^*})}_G \,\mm{ and } \, \widehat{\Sigma} \cong\p_G(H \setminus \Sigma,d_{\bp}) \cong \p_G(H \setminus \Sigma,d_{\bp^*}),
  \]
  where for $X=(H \setminus \Sigma,d_{\bp})$ or $(H \setminus \Sigma,d_{\bp^*})$, $\overline{X}_G$ and $\p_G(X)$ denote the Gromov compactification and the Gromov boundary respectively. For $\widehat{H}$ and $\widehat{\Sigma}$, see Section~\ref{notation} after Remark~\ref{ram}.
\end{itemize}
The spaces $(H \setminus \Sigma, d_{\bp})$ and $(H \setminus \Sigma, d_{\bp^*})$ are conformally equivalent to the original space $(H \setminus \Sigma, g_H)$. We refer to both these spaces as \textbf{hyperbolic unfoldings} of $(H \setminus \Sigma, g_H)$.
\end{theorem}

The completeness and the bounded geometry property are consequences of the \si-axioms. The Gromov hyperbolicity of $(H\setminus \Sigma, d_{\bp})$ relies on the \si-uniformity of $H \setminus \Sigma$; ordinary uniformity would not be sufficient to guarantee the hyperbolicity of the \si-metric. On the other hand, the classical \emph{quasi-hyperbolic metric} on $H \setminus \Sigma$ (cf.~\cite{BHK}) usually has non-bounded geometry and will not be used in our papers cf. the discussion before \cite[Theorem 1.11]{L1}.

\begin{remark}[Quantitative Control]\label{qic}
The geometry of the hyperbolic unfoldings of $H \in {\cal{G}}$ described in \ref{thm3} admits some quantitative estimates. We summarize results from \cite[2.7(iv), 3.6, 3.7 and 3.11]{L1} and recall that the \si-uniformity constant
 of $(H \setminus \Sigma, g_H)$
largely determines the hyperbolicity constant of $(H \setminus \Sigma, d_{\bp})$. For the notationally simplest case where the Lipschitz constant $L_{\bp}$ of $\delta_{\bp}$ equals $1$, which can be realized as $\bp_1$ in  \cite[Definition 2.2]{L1}, we have
\begin{itemize}
  \item When $H$ is $a$-\si-uniform, the \si-metric $d_{\bp}$ is $\delta(a)$-hyperbolic with
  \begin{equation}\label{du}
  \delta(a)=4a^2 \cdot \log \big(1+  c(a) \cdot ( 2c(a) + 3)\big).
  \end{equation}
  \item The Whitney smoothed \si-metric $d_{\bp^*}$ is $\Delta$-hyperbolic with $\Delta(H)$.
\end{itemize}
For $H \in {\cal{H}}^{\R}_n$,  we even have that $H$ is  $c_n$-\si-uniform for a  constant $c_n>0$ depending only on the dimension $n$ and we get constants $\delta_n ,\Delta_n >0$ depending only on $n$ so that
\begin{itemize}
  \item   $(H \setminus \Sigma,d_{\bp})$ is $\delta_n$-hyperbolic and  $(H \setminus \Sigma,d_{\bp^*})$ is $\Delta_n$-hyperbolic.
\end{itemize}
Moreover, for $H \in {\cal{H}}^{\R}_n$, we get bounds depending only on $n$ for the constants, in  Definition~\ref{bog}, which quantify the boundedness of the geometry of $(H \setminus \Sigma, d_{\bp})$ and $(H \setminus \Sigma, d_{\bp^*})$:
\begin{itemize}
  \item $\varrho_{(H \setminus \Sigma, d_{\bp})} = a_ {\R}(n)$, $\ell_{(H \setminus \Sigma, d_{\bp})} = b_ {\R}(n)$   and  $\varrho_{(H \setminus \Sigma, d_{\bp^*})} = a^*_ {\R}(n)$,  $\ell_{(H \setminus \Sigma, d_{\bp^*})} = b_ {\R}^*(n)$.
\end{itemize}
Finally, note that $d_{\bp}$ does not change under scalings of the underlying space $H \in {\cal{G}}$ by some constant $\tau >0$, since the transformed arc length will be multiplied by $\tau$ and $\bp$ divided by $\tau$. In particular, for any given $H \in {\cal{G}}$ we get a common bound on $\varrho$ and $\ell$ for $\tau \cdot H$, $\tau \ge 1$,  and its blow-up geometries like tangent cones.  \qed
\end{remark}

%
%
\subsubsection{Hyperbolic Boundary Harnack Principle}\label{app}

Next we  describe the range of admissible elliptic operators:

\begin{definition}\label{sao0}
For a complete Riemannian manifold $X$ with bounded geometry, we call a second order elliptic operator $L$ on $X$ \textbf{adapted weakly coercive} provided the following conditions hold:
\begin{itemize}
  \item $L$ is \textbf{adapted}: There exists a constant $k=k_L\geq 1$ such that $L$ satisfies the following uniformity conditions. With respect to the charts $\phi_{p}$, $L$ can be locally written as
  {\small \[
  -L(u) = \sum_{i,j} a_{ij} \cdot \frac{\p^2 u}{\p x_i \p x_j} + \sum_i b_i \cdot \frac{\p u}{\p x_i} + c \cdot u
  \]}
  for $\beta$-H\"{o}lder continuous $a_{ij}$, $\beta\in (0,1]$, and measurable functions $b_i$, $c$, such that for all $p\in X$, $\xi\in\R^n$, we have
  \begin{enumerate}
    \item $k^{-1} \cdot\sum_i \xi_i^2 \le \sum_{i,j} a_{ij}(p) \cdot \xi_i \xi_j\le k \cdot\sum_i \xi_i^2$,
    \item $|a_{ij}|_{C^\beta(B_{\rho}(p))}\le k$,  $|b_i|_{L^\infty}$ and $|c|_{L^\infty} \le k.$
  \end{enumerate}
  \item $L$ is \textbf{weakly coercive}: There exists a positive $C^2$-supersolution $u$ of the equation $L \, f = 0$ with $L \, u \ge \ve \cdot u$ for some $\ve >0$.
\end{itemize}
\end{definition}

One readily checks that there is a largest such $\ve>0$, the \textbf{(generalized) principal eigenvalue} $\tau_L=\tau(L,X)>0$. Moreover, there is a smallest such $k$ written  $\kappa_L=\kappa(L,X)\ge 1$.\\

We notice that both $L$ and $L - \ve \cdot Id$ are adapted weakly coercive if $\ve <\tau_L$, whereas $L - \tau_L \cdot Id$ is no longer weakly coercive. \\

For adapted weakly coercive operators on Gromov hyperbolic manifolds of bounded geometry Ancona has proved a \emph{hyperbolic} BHP relative to $\p_GX$. The boundary condition corresponding to the vanishing of $u$ and $v$ in the \emph{euclidean} BHP (cf.\ the beginning of Section~\ref{pre}) is that of \emph{L-vanishing}: A solution $u \ge 0$ \textbf{\emph{L}-vanishes} along some open subset $V$ of $\p_GX$ if there is a supersolution $w >0$  such that $u/w \ra 0$ when we approach $V$ in $X$.\\

Canonical $\Phi_\delta$-chains are well-suited for the following central result of Ancona's potential theory on Gromov hyperbolic manifolds of bounded geometry  cf.\cite[Cor.5.8]{KL}.
\begin{theorem}[Ancona's Hyperbolic BHP]\label{abxhp}
Let $X$ be a complete $\delta$-hyperbolic manifold of $(\varrho,\ell)$-bounded geometry and assume that $L$ is an adapted weakly coercive operator on $X$.\\

 Let  $u$, $v>0$ be two supersolutions* of $L\, f = 0$ on $X$ both properly solving $L\, f = 0$ on $\mathcal{N}^\delta_i $ and $L$-vanishing along $\mathbf{N}^\delta_i \cap \p_GX$, then there exists a constant $C=C(\varrho,\ell,k,\tau,n,\delta)>1$ solely depending on $\varrho$, $\ell$, $k$, $\tau$, $n$ and $\delta$ such that
\begin{equation}\label{fhepq}
u(x)/v(x) \le C \cdot  u(y)/v(y), \mm{ for all } x,\, y \in \mathcal{N}^\delta_{i+1}  \mm{ and any $i =1,....m-1$}.
\end{equation}
\end{theorem}

 *A supersolution is a lower semi-continuous function with values in $\R \cup \{+\infty\}$ that is larger than the Dirichlet solution on any ball with the same boundary values and finite on a dense set.
 An important example is the minimal Green's function $G(\cdot,p)$ of $L$. It is a supersolution $L$-vanishing along the entire boundary $\p_GX$.

\begin{remark} \label{regg} This version of the BHP follows from a step-by-step review of Ancona's arguments in \cite{An1} and \cite{An2} with some refinements taken from \cite{BHK}. The details are presented in \cite{KL}. In  ~\cite{KL} we impose a stronger regularity on the coefficients than in Def.\ref{sao0} to simplify the exposition. We use $a_{ij} \in C^{2,\beta}, b_i \in C^{1,\beta}, c \in C^{\beta}$. (These conditions are sufficient for geometric and physical applications, as in ~\cite[Theorem 2]{L2}.)  However,\ref{abxhp} still holds under the weaker conditions in Def.\ref{sao0} with technical adjustments explained in \cite[Ch.5]{An1} and \cite{An2}. \qed
\end{remark}

%
%
%
\setcounter{section}{3}
\renewcommand{\thesubsection}{\thesection}
\subsection{Boundary Harnack Principles and Martin Theory on $H \setminus \Sigma$}
%

We first turn to the potential theory of \si-adapted operators on $H \setminus \Sigma$. As indicated in the introduction we use the \si-uniformity of $H \setminus \Sigma$ to unfold almost minimizers to complete hyperbolic spaces. Here, we have boundary Harnack inequalities for uniformly elliptic and even more general operators.

%
\subsubsection{\si-Adapted Operators on $(H \setminus \Sigma, g_H)$}\label{agrm}
The results in the last section show that any adapted weakly coercive operator on the unfolding $X=(H \setminus \Sigma, d_{\bp^*})$ satisfies Ancona's BHP. If we translate this back to the original space $(H \setminus \Sigma, g_H)$ we arrive at the notion of \si-adapted operators.\\

\noindent\textbf{\si-Adapted Operators.} \,
Given $H \in \cal{G}$, we will work with an \si-adapted atlas of $H$ consisting of \textbf{\si-adapted charts}: For any $K>1$, there is a radius $\Theta(p):= \Gamma / \bp(p)$ with $\Gamma(H,K,\bp) >0$ such that for any $p \in H \setminus \Sigma$, the exponential map $\exp_p|_{B_{\theta(p)}(0) \subset T_pH}$ is a $K$-bi-Lipschitz
$C^\infty$-diffeomorphism onto its image, cf.\ \cite[Proposition B.1]{L1}. Thus, we get the smooth charts
\[
\psi_{p}:=\exp^{-1}_p|_{B_{\Theta(p)}(0)}: B_{\Theta(p)}(p) \ra B_{\Theta(p)}(0) \subset \R^n, \quad \psi_{p}(p)=0
\]
on $H\setminus\Sigma$. As a direct counterpart or reformulation of \ref{qic}  we also get a quantitative control over $\Gamma$. In particular, when $H \in {\cal{H}}^{\R}_n$  we have, for any fixed $\bp$, $\Gamma=\Gamma(K,n)>0$.
In what follows the default choices are $K=2$ and $\Gamma(H,2,\bp)$.

\begin{definition}\label{sao2}
Let $H \in \cal{G}$. A second order elliptic operator $L$ on $H\setminus \Sigma$ is called \textbf{\si-adapted} if the following two conditions hold.
\begin{itemize}
  \item \emph{\textbf{$\mathbf{\bp}$-adaptedness:}} \,
  There exists a constant $k=k_L\ge 1$ such that $L$ satisfies the following \si-weighted uniformity conditions. With respect to the charts $\psi_p$, $L$ can be locally written
  {\small \[
  -L(u) = \sum_{i,j}  a_{ij} \cdot \frac{\p^2 u}{\p x_i \p x_j} + \sum_i b_i \cdot \frac{\p u}{\p x_i} + c \cdot u
  \]}
  for $\beta$-H\"{o}lder continuous coefficients $a_{ij}$, $\beta\in (0,1]$, and measurable functions $b_i$, $c$, such that for all $p\in H\setminus\Sigma$, $\xi\in\R^n$, we have
  \begin{enumerate}
    \item $k^{-1} \cdot\sum_i \xi_i^2 \le \sum_{i,j} a_{ij}(p) \cdot \xi_i \xi_j \le k \cdot \sum_i \xi_i^2$,
    \item $\delta^{\beta}_{\bp}(p) \cdot  |a_{ij}|_{C^\beta(B_{\theta(p)}(p))} \le k$, $\delta_{\bp}(p) \cdot |b_i|_{L^\infty} \le k$ and $\delta^2_{\bp}(p) \cdot |c|_{L^\infty} \le k$.
  \end{enumerate}
  \item \emph{\textbf{$\mathbf{\bp}$-weak coercivity:} \,
  There exists a positive $C^2$-supersolution $u$ of the equation $L \, f = 0$ with}
  \[
  L \, u \ge \ve \cdot \bp^2 \cdot u\mm{ for  some } \ve  >0.
  \]
\end{itemize}
\end{definition}

As in \ref{sao0} there is a largest such $\ve>0$, the  principal eigenvalue  $\tau_L=\tau(L,(H \setminus \Sigma, g_H))>0$. Moreover, there is a smallest such $k$ written  $\kappa_L=\kappa(L,(H \setminus \Sigma, g_H))\ge 1$.

\begin{remark}
In \cite{L2} we will show that \si-adapted operators exist in abundance by using a so-called \emph{Hardy structure} on $H \in {\cal{G}}$.\qed
\end{remark}

First we show that the notion of \si-adaptedness on $(H \setminus \Sigma, g_H)$ is the counterpart of adapted weak coercivity on $(H \setminus \Sigma, d_{\bp^*})$.

\begin{proposition}[Unfolding Correspondence]\label{ucor}
For any $H \in {\cal{G}}$, we consider the canonical correspondence
\[
(H \setminus \Sigma, g_H)  \mm{ equipped with }
L \, \rightleftharpoons \, (H \setminus \Sigma, d_{\bp^*})\mm{ equipped with } \delta_{\bp^*}^2 \cdot L,
\]
Then we have the equivalence
\begin{equation}\label{unfold}
L \mm{ is \si-adapted} \,\, \Leftrightarrow \,\, \delta_{\bp^*}^2\cdot L \mm{ is adapted weakly coercive}.
\end{equation}
For the two constants $\tau$ and $\kappa$ we have:
\begin{equation}\label{unfold1}
c_*^{-1} \cdot \tau(L,(H \setminus \Sigma, g_H)) \le \tau(L,(H \setminus \Sigma, d_{\bp^*}) \le c_* \cdot \tau(L,(H \setminus \Sigma, g_H)),
\end{equation}
\begin{equation}\label{unfold2}
c_*^{-1} \cdot \kappa(L,(H \setminus \Sigma, g_H)) \le \kappa(L,(H \setminus \Sigma, d_{\bp^*}) \le c_* \cdot \kappa(L,(H \setminus \Sigma, g_H))
\end{equation}
for some constant $c_*(L_{\bp},H) \ge 1$, in general,  and for $H \in {\cal{H}}^{\R}_n$ we have $c_*(L_{\bp},n) \ge 1$.
\end{proposition}
\noindent\textbf{Proof} \,
Assume that $L$ satisfies the \si-adaptedness conditions in \ref{sao2}. Then we consider the \emph{locally constantly} scaled version of $\psi_{p}$
\[
\psi^{\bp}_{p}:(B_{\gamma}(p),\bp^2(p) \cdot g_H) \ra \R^n\quad\mm{with}\quad
\psi^{\bp}_{p}(x):= \bp(p) \cdot \psi_{p}(x).
\]
Since both the source and the target have been scaled by the same constant, $\psi^{\bp}_{p}$ is again a $2$-bi-Lipschitz map.
These particular scalings given rise to a neat transformation law for the given \si-adapted operator $L$. Indeed, if we denote by $y_i$ the coordinates induced by
$\psi^{\bp}$, then $-L \, u$ can be recomputed, using the chain rule, as
\[
\bp^2(p) \cdot \sum_{i,j} a_{ij} \cdot \frac{\p^2 u}{\p y_i \p y_j} + \bp(p) \cdot \sum_i b_i \cdot \frac{\p u}{\p y_i} + c \cdot u.
\]
The adaptedness of $L$ with respect to $\bp$ and the Lipschitz continuity of both $\psi^{\bp}_{p}$ and $\delta_{\bp}$ then show that $\delta_{\bp^*}^2\cdot L$ satisfies the adaptedness condition
for adapted weakly coercive operators on $(H \setminus \Sigma, d_{\bp^*})$ with respect to the charts $\psi^{\bp}_{p}$. Next, the weak coercivity of $L$ relative to $\bp$ implies that there exists a positive supersolution $u$ of the equation $L \, f = 0$. Hence,
\[
\delta_{\bp^*}^2\cdot L \, u \ge \ve \cdot
\delta_{\bp^*}^2\cdot \bp^2 \cdot u \mm{ for  some } \ve  >0.
\]
Thus the approximation property $c_1 \cdot \delta_{\bp}(x) \le \delta_{\bp^*}(x)  \le c_2 \cdot
\delta_{\bp}(x)$,for the Whitney smoothing constants $c_1, c_2$ in \ref{smsk}, gives the weak coercivity of $\delta_{\bp^*}^2\cdot L$ on $(H \setminus \Sigma, d_{\bp^*})$. \\

Finally, the two
characteristic constants $\tau$ and $\kappa$ of $L$ relative $(H \setminus \Sigma, g_H)$ and of $\delta_{\bp^*}^2\cdot L$  relative $(H \setminus \Sigma, d_{\bp^*})$  differ by multiples $c_* \ge 1$
depending only on $c_1, c_2$ and the Lipschitz constant $L_{\bp}$. We have the dependencies $c_*=c_*(L_{\bp},H)$, in general, and $c_*=c_*(L_{\bp},n)$, for $H \in {\cal{H}}^{\R}_n$.\\

For the converse we can argue in the same way and pass from adapted weakly coercive to \si-adapted operators.\qed

The hyperbolic unfoldings combined with the hyperbolic BHP  gives us the corresponding BHP for \si-adapted operators directly on the original $H \in {\cal{G}}$ where
$\Sigma_H$ is considered as the boundary. We can now (re)formulate the versions of \ref{abxhp} on the original almost minimizer we get from the unfolding correspondence.

\begin{theorem}[BHP for \si-Adapted Operators]\label{mbhsq}
Let $H \in {\cal{G}}$ and $L$ be an \si-adapted operator on $H \setminus \Sigma$.
Let  $u$, $v>0$ be two supersolutions of $L\, f = 0$ on $H \setminus \Sigma$ both solving $L\, f = 0$ on $\mathbf{N}^\delta_i \cap H \setminus \Sigma$ and $L$-vanishing along $\mathbf{N}^\delta_i \cap \widehat{\Sigma}$, then there is a $C=C(H,L)>1$ such that
\begin{equation}\label{fhepq}
u(x)/v(x) \le C \cdot  u(y)/v(y), \mm{ for all } x,\, y \in \mathcal{N}^\delta_{i+1} \mm{ and any $i =1,....m-1$},
\end{equation}
where $\delta$ is the hyperbolicity constant of the hyperbolic unfolding from \ref{qic}.
\end{theorem}

On the geometric side, \ref{abxhp} shows that the constant $C(H,L)$ only depends on the hyperbolicity and the bounded geometry constants and the dimension of its hyperbolic unfolding.
Remarkably, the hyperbolic unfoldings of all Euclidean hypersurfaces $H \in {\cal{H}}^{\R}_n$ satisfy the same estimates for the hyperbolicity and the bounded geometry constants cf.\ref{qic}. Therefore we get the following refinement of \ref{mbhsq}.

\begin{theorem}[Stable BHP for \si-Adapted Operators]\label{mbhsqq} Under the assumptions of \ref{mbhsq} we have for $H \in {\cal{H}}^{\R}_n$ and $L$ with $k_L \le \kappa$ and $\ve_L \ge \eta$ the constant $C$ only depends on $\kappa$, $\eta$ and $n$.\\
\end{theorem}
%
%
\subsubsection{Martin Theory for \si-adapted Operators}\label{nmt}
Let $X$ be a non-compact Riemannian manifold and $L$ be a linear second order elliptic operator on $X$. The local BHP links the geometry of the boundary (its uniformity) to representations of minimal solutions of $L\,f=0$ in terms of Martin (boundary) integrals. We start with a short reminder of some basics from Martin theory, cf.\ \cite[Chapter I.7]{BJ},  \cite[Section 7.1]{P} for details.\\

Recall first that a Green's function $G$ is \textbf{minimal} if there is no distinct positive solution $w$ on $X$ with $w \le G(\cdot,y)$. This is tantamount to saying that $G(\cdot,x)$ is an $L$\textbf{-potential}.

\begin{definition}[Martin Boundary]\label{mb}
Let $X$ be a non-compact Riemannian manifold and $L$ be a linear second order elliptic operator on $X$ with a minimal Green's function $G: X\times X \ra (0,\infty]$. We choose a base point $p$ and consider the space $S$ of sequences $s=\{p_n\}$ in $X$, $n \ge 1$, such that
\begin{itemize}
  \item $s$ has no accumulation points in $X$.
  \item $K(x,p_n):= G(x,p_n)/G(p,p_n) \ra K_s(x)$ compactly to some function $K_s$ on $X$ as $n \ra \infty$.
\end{itemize}
The \textbf{Martin boundary} $\p_M (X,L)$ is the quotient of $S$ modulo the following relation on $S$: $s \sim s^*$ if and only if $K_s \equiv K_{s^*}$. Moreover, we define the \textbf{Martin kernel} $k(x;y)$ on $X \times \p_M (X,L)$ by $k(x;y):= K_s(x)$, for some sequence $s$ representing $y \in \p_M (X,L)$. As for the Gromov boundary, these definitions do not depend on the choice of the base point $p$.
\end{definition}

\begin{remark}\label{mrtt}
The Harnack inequality and elliptic theory show that each $K_s \in \p_M (X,L)$ is a positive solution of $L \, u = 0$ on $X$. Thus the convex set
\[
S_L(X)=\{u\in C^{2,\beta}(X)\mid L\, u =0,\,u>0,\,u(p) =1\}
\]
is compact in the topology of compact convergence, and $\p_M (X,L)$ is a compact subset of $S_L(X)$.\qed
\end{remark}

The (metrizable) \textbf{Martin topology} on $\overline{X}_M:= X \cup \p_M (X,L)$ is defined as follows (cf.\ \cite[Chapter I.7]{BJ} or \cite[Chapter 12]{H} for more details). A sequence $s=\{p_n\}$ with no accumulation points in $X$ converges to a point $y \in \p_M (X,L)$ if and only if $K(x,p_n)$ converges compactly to $K_s(x)$ representing $y$. Furthermore, $y_n \in  \p_M (X,L)$ converges to $y  \in \p_M$ if and only if the functions $K_s(y_n)(x)$ representing $y_n$ converge compactly to functions $K_s(y)(x)$ representing $y$. It follows that $\overline{X}_M$ is compact, and $\p_M (X,L)$ is a closed subspace. The space $\overline{X}_M$ is also called the \textbf{Martin compactification} of $(X,L)$.\\

\textbf{Minimal Martin boundary.} A solution $u\in S_L(X)$ of $L \, f = 0$, is called extremal if it cannot be written as a non-trivial convex combination of other elements in $S_L(X)$. It is easy to see that $u$  is extremal  if and only if $u$ is a minimal solution.  One calls $u$ \textbf{minimal} if for any other solution $v >0$ with $v \le u$, we have $v \equiv c \cdot u$ for some constant $c>0$. The subset $\p^0_M (X,L) \subset \p_M (X,L)$ of extremal points of $S_L(X)$ is called the \textbf{minimal Martin boundary}. Typically, it is much smaller than the full Martin boundary.\\

Since the $\mathbf{N}^\delta_i$  form a neighborhood basis of each singular point, the BHP   implies, from  standard comparison arguments that the (minimal) Martin boundary $\p^0_M (X,L)$ is the topological respectively ideal boundary $\p X$ of $X$. Details are explained, for instance, in \cite{KL},  \cite{Ai1}, \cite{AG} or \cite{An1}.

\begin{theorem}[Martin Boundary for \si-Adapted Operators]\label{mbhs}
Let $H \in {\cal{G}}$ and $L$ be an \si-adapted operator on $H \setminus \Sigma$. Then
\begin{itemize}
  \item the identity map on $H \setminus \Sigma$ extends to a homeomorphism between $\widehat{H}$ and the Martin compactification $\overline{(H \setminus\Sigma)}_M$.
  \item all Martin boundary points are minimal:  $\p^0_M (H \setminus \Sigma,L) \equiv \p_M(H \setminus \Sigma,L)$.
\end{itemize} Thus,  $\widehat{\Sigma}$ and the minimal Martin boundary $\p^0_M (H \setminus \Sigma,L)$ are homeomorphic. In the case where $H$ is a singular cone $H=C \in \mathcal{SC}_{n}$, this, in particular, means that there is exactly one (minimal) Martin boundary point at the origin and one at infinity.
\end{theorem}

\begin{remark}\label{alt}
Alternatively, we may prove Proposition~\ref{mbhs} by starting from the Martin theory of the hyperbolic unfolding $X=(H \setminus \Sigma, d_{\bp^*})$ and using the hyperbolic BHP. Indeed, the identity map $id_{H \setminus \Sigma}$ extends to homeomorphisms between $H$ resp.\ $\widehat{H}$, the Gromov compactification $\overline{X}_G$, and the Martin compactification $\overline{X}_M$ of $X$ for any adapted weakly coercive operator, in our case $\delta_{\bp^*}^2 \cdot L$. See \cite[Theorems 2 and 8]{An1}, \cite[Theorem V.6.2 p.\ 97]{An2} and also \cite[Chapter 8]{BHK} for details. In a nutshell, these results translate from $(H \setminus \Sigma, d_{\bp^*})$ to $(H \setminus \Sigma, g_H)$ by applying the unfolding correspondence \eqref{unfold} as the positive factor  $\delta_{\bp^*}^2$ does neither change the notion of positive nor minimal solution for $L \, f = 0$. Put differently, up to obvious adjustments in the statements, the respective BHPs for positive solutions of $L \, f = 0$ as well as the respective Martin compactifications and boundaries are simply identical, that is, $\overline{(H \setminus \Sigma, g_H)}_M \equiv \overline{(H \setminus \Sigma, d_{\bp^*})}_M$
In particular, we have $\p_M((H \setminus \Sigma, g_H),L)  \equiv \p_M((H \setminus \Sigma, d_{\bp^*}),\delta_{\bp^*}^2 \cdot L)$.\qed
\end{remark}

\begin{remark}[Almost Minimizers with Boundary]\label{mgr}
We get similar results for area minimizers $H$ with boundary $\p H$, that is, $H$ solves the \textbf{Plateau problem} for the boundary $\p H$ with $\Sigma \cap \p H \v$. For this we replace $\delta_{\bp}(x)$ by $\mathbf{d}(x):=\min \{dist(x,\p H),\delta_{\bp}(x)\}$. Subject to
\[
H\mm{ is a \emph{uniform space} \,  and  \, } dist(z,\Sigma) \le c_H \cdot \delta_{\bp}(z)
\]
for some constant $c_H>0$ one deduces a version of \si-uniformity as in Theorem~\ref{thm2} after replacing $\delta_{\bp}$ by $\mathbf{d}$. Then we construct hyperbolic unfoldings for $\mathbf{d}$ (and similarly for its Whitney smoothing $\mathbf{d}^*$) for the distance function
\[
d_{\mathbf{d}}(x,y) := \inf \Bigl  \{\int_\gamma  1/\mathbf{d}(\cdot) \, \, \Big| \, \gamma   \subset  H \setminus \Sigma\mbox{ rectifiable curve joining }  x \mbox{ and } y  \Bigr \}.
\]
As in \cite{L1} and \cite{BHK} we see that $d_{\mathbf{d}}$ and $d_{\mathbf{d}^*}$ define complete Gromov hyperbolic spaces with bounded geometry such that
\[
\widehat{H} \cong \overline{(H \setminus \Sigma,d_{\mathbf{d}})}_G  \cong \overline{(H \setminus \Sigma,d_{\mathbf{d}^*})}_G\quad\mm{and}\quad\widehat{\Sigma \cup \p H}\cong\p_G(H \setminus \Sigma,d_{\mathbf{d}}) \cong \p_G(H \setminus \Sigma,d_{\mathbf{d}^*}).
\]
We have now \emph{two} boundary components in the Gromov boundary corresponding to $\Sigma$ and $\p H$ respectively. On the analytic side we consider a second order elliptic operator $L$ on $D$ with the appropriate adaptedness properties:
\begin{itemize}
  \item $L$ satisfies $\mathbf{d}$-weighted uniformity conditions: Locally,
  {\small \begin{equation}\label{ada}
  -L(u) = \sum_{i,j} a_{ij} \cdot \frac{\p^2 u}{\p x_i \p x_j} + \sum_i b_i \cdot \frac{\p u}{\p x_i} + c \cdot u
  \end{equation}}
  with $\beta$-H\"{o}lder continuous coefficients $a_{ij}$, $\beta\in(0,1]$, and measurable functions $b_i$, $c$. Further, there exists a $k \ge 1$ such that for any $p\in H\setminus\Sigma$, $\xi\in\R^n$, we have
  \begin{enumerate}
    \item $k^{-1} \cdot\sum_i \xi_i^2\le \sum_{i,j} a_{ij}(p) \cdot \xi_i \xi_j \le k \cdot \sum_i \xi_i^2$,
    \item $\mathbf{d}(p)^{\beta} \cdot  |a_{ij}|_{C^\beta(B_{\theta(p)}(p))} \le k$, $\mathbf{d}(p) \cdot |b_i|_{L^\infty} \le k \mm{ and } \mathbf{d}(p)^2 \cdot |c|_{L^\infty} \le k$.
  \end{enumerate}
  \item There exists a positive $C^2$-supersolution $u$ of the equation $L \, f = 0$ with
  \[
  L \, u \ge \ve \cdot \mathbf{d}^{-2} \cdot u, \mm{ for  some } \ve  >0.
  \]
  Then, as in Proposition~\ref{mbhs}, the identity map on $H$ extends to a homeomorphism between $\widehat{H}$ and the Martin compactification, and all Martin boundary points are minimal:
  \[
  \widehat{H} \cong \overline{H \setminus \Sigma}_M \quad  \mm{ and } \quad \widehat{\Sigma \cup \p H}\cong \p^0_M (H\setminus\Sigma,L).
  \]
\end{itemize}
\noindent One gets similar extensions also for the other results of this paper.\qed
\end{remark}

\setcounter{section}{4}
\renewcommand{\thesubsection}{\thesection}
\subsection{Basic Asymptotic Analysis}
The boundary Harnack principles and the resulting Martin Theory on $H \setminus \Sigma$ have deep consequences for the asymptotic analysis of \si-adapted operators. In this chapter we consider some classical boundary value problems like Fatou's theorem, or the Dirichlet problem on $H \setminus \Sigma$ regarding $\Sigma$ as a boundary. The discussion involves new ideas, for instance the \emph{freezing of tangent cone approximations of almost minimizers}. This will be also an essential input in the follow-up paper~\cite{L2}.

\subsubsection{Tangent Cones and Freezings}\label{tfr}
Almost minimizers can be approximated by area minimizing cones around singular points. We review these approximations and use \si-structures to derive some refinements building on the \emph{ asymptotic freezing effect}.\\

\noindent\textbf{Spaces of Tangent Cones.} \,
Let $H \in {\cal{G}}$. Blowing up at $p\in \Sigma_H$, that is, rescaling $H$ around $p$ by a strictly increasing sequence $\tau_m\to\infty$, yields a subconverging sequence $\tau_{m_k}\cdot H$ whose limit is an again area minimizing cone, a so-called \emph{tangent cone}. It can be viewed as a partially linearized and simplified local model for $H$. For a formal definition, cf.\ \cite[4.3.16]{F}, \cite[Chapter 37.4]{Si}, \cite[Theorem 1]{T1} and also \cite[Appendix A]{L1} for a short review of related results such as the following

\begin{proposition}[Convergence to Tangent Cones]\label{flat-norm-approx}
Let $H \in {\cal{G}}$ and $p \in \Sigma_H$. For every sequence $\tau_m \to +\infty$ of positive real numbers there exists a subsequence $\tau_{m_k}$, as well as an area minimizing cone $C_p \subset \R^{n+1}$, with $0 \in \sigma_C$, such that
\begin{itemize}
  \item \textbf{flat norm convergence:} For any given open $U \subset \R^{n+1}$ with compact closure the \textbf{flat norm} $\db_U$ converges to zero: \, $\db_U (\tau_{m_k} \cdot H, C_p) \to 0.$
  \item  \textbf{$C^l$-norm convergence:} If, in addition, $\overline U \subset C_p \setminus \sigma$, then $\db_U$-convergence implies compact $C^l$-convergence, for any $l \ge 0$.
\end{itemize}
Recall that $\sigma$ is our generic notation for singularity sets of cones, cf.\ the introduction. The (pseudo-)metric $\db_U$ can be thought of as the volume between the two hypersurfaces in $U$. For $H \in {\cal{G}}$,   $\db_U (\tau_{m_k} \cdot H, C_p) \to 0$ is equivalent to the condition that the Hausdorff-distance $d_{GH}(U \cap (\tau_{m_k} \cdot H), U\cap  C_p) \ra 0$.
\end{proposition}


\begin{remark}[Iterated Tangent Cones]\label{itc}
When we blow up around a point $p_0\in \Sigma_H$ we get a singular tangent cone $C$. Scaling around the tip $0 \in \sigma_{C}$ merely reproduces $C$, but blowing up singular points $p_1\neq 0 \in \sigma_{C_1}$ gives rise to \emph{iterated tangent cones}. The resulting area minimizing cones $C'$ are Riemannian products $\R \times C^{n-1}$, where $C^{n-1}$ is an area minimizing cone in $\R^n$. This iteration process ends when we arrive at a cone $\R^m \times C^{n-m} $ where $C^{n-m} \subset \R^{n-m+1}$ is singular only at $0$. It is a basic fact that $n-m \ge 7$, since minimal cones of dimension equal or less than six are necessarily regular, cf.\ \cite[Chapter 11]{Gi}.\qed
\end{remark}

Let ${\cal T}_p \subset\mathcal{SC}_n$ be the set of  tangent cones at $p \in \Sigma \subset H$ and ${\cal T}_H=\bigcup_{p\in\Sigma}{\cal T}_p$ be the set of tangent cones of points in $\Sigma$. Note that together with a tangent cone $C$ we include its image under any linear isometry of $\R^{n+1}$, that is, the subsets ${\cal T}_p$ and ${\cal T}_H$ of $\mathcal{SC}_n$ are $O(n+1)$-homogenized.\\

A major problem is the occurence of possibly infinitely many distinct tangent cones at $p \in \Sigma$ as different subsequences of $\tau_{m} \cdot H,$ may converge to distinct limits. Also, ${\cal T}_p$ and the control on the sequence $\db_U (\tau_{m_k} \cdot H, C_p)$ may change discontinuously in $p$. Fortunately, we can invoke the following compactness results for ${\cal T}_p$ and ${\cal T}_H$ which are obtained from Propostion~\ref{flat-norm-approx} and the definition of tangent cones.

\begin{lemma}[Compact Cone Spaces]\label{cone-compactness}
Under compact convergence in flat norm topology, the following statements are true.
\begin{enumerate}
  \item The spaces $\mathcal{C}_{n}$ and ${\cal{H}}^{\R}_n$ are compact. Further, $\mathcal{K}_{n}$ is compact in flat norm topology.
  \item There exists a constant $d_n >0$ such that for $C \in \mathcal{C}_{n}$ we have
  \[
  \db_{B_1(0)} (C, Y) < d_n\mm{ for some hyperplane } Y \subset \R^{n+1} \Leftrightarrow C \mm{ is non-singular.}
  \]
  \item The sets $\mathcal{SC}_n \subset\mathcal{C}_n$ and ${\cal{SH}}^{\R}_n \subset {\cal{H}}^{\R}_n$ are closed and thus compact. Consequently, $\overline {\cal T}_H \subset \mathcal{SC}_n$ and $\overline {\cal T}_p =  {\cal T}_p$.
\end{enumerate}
\end{lemma}

\begin{remark}[Notions of Compactness]\label{noc}
As in the previous lemma the general compactness results for currents refer to \emph{sequential compactness} relative to the topology given by the family of flat pseudometrics. That is, $H_i\to H$ if and only if $\db_U(H_i,H)\to0$ for all bounded open $U\subset\R^{n+1}$, cf.\ \cite[Chapter 31]{Si}. Note, however, that each $N\in \mathcal{K}_{n}$ has a unique representation as the compact metric completion of the manifold $N \setminus \Sigma_N$. Thus, on $\mathcal{K}_{n}$, we have a proper flat \emph{metric}. Here, sequential compactness therefore implies compactness, i.e., each open cover of $\mathcal{K}_{n}$ contains a finite subcover.\qed
\end{remark}

\noindent\textbf{\si-Freezing Effects near Singular Points} \,
Now we show that tangent cones can still serve as local models for area minimizers near a singular point $p$. For this, we measure the \emph{non-tangential approachability} of $p$.

\begin{definition}\label{pen}
For $H \in {\cal{G}}$ and  $\omega>0$, we define the \textbf{\si-pencil $\P(p,\omega)$ pointing to $p\in \Sigma$} by
\begin{equation}\label{ii}
\P(p,\omega)=\P_H(p,\omega):= \{x \in H \setminus \Sigma \,|\, \delta_{\bp} (x) > \omega \cdot d_H(x,p)\}.
\end{equation}
\end{definition}

Note in passing the scaling invariance $\P_{\tau\cdot H}=\P_H$ which follows from \si-transform axiom (S2). In view of scaling arguments it is also useful to consider the \textbf{truncated \si-pencil}
\[
\TP(p,\omega,R,r)=\TP_H(p,\omega,R,r):=B_{R}(p) \setminus B_{r} (p) \cap  \P(p,\omega)\subset H.
\]
While we zoom into some singular point by rescaling $\tau \cdot  \TP_H(p,\omega,R/\tau,r/\tau)$, the \si-pencil $\P(p,\omega)$ is better and better $C^k$-approximated by (usually changing) tangent cones while the twisting of $\P(p,\omega)$ slows down as $\tau\to\infty$.

\begin{proposition}[Asymptotic \si-Freezing]\label{freez}
Let $H\in\cal{G}$ and $p \in \Sigma_H$. Further, pick $\ve > 0$ and a pair $R > 1 > r >0$. Then the following is true.
\begin{itemize}
  \item \emph{\textbf{Flat-norm Version:}} \,
  We have some $\tau_{\ve , R , r,p} > 1$ such that for every $\tau \ge \tau_{\ve , R , r,p}$ there is a tangent cone $C_p^\tau \in \mathcal{SC}_n$ of $H$ at $p$ with
  \[
  \tau \cdot (H \cap B_{R/\tau}(p) \setminus B_{r/\tau} (p)) \text{ is $\ve$-close in flat norm to } C^\tau_p\cap B_R(0) \setminus B_r (0 ).
  \]
  \item \emph{\textbf{$\mathbf{C^k}$-norm Version:}} \,
  For any additionally given $1 > \omega >0$ and $k \in \Z^{\ge 0}$ we can find $\tau_{\ve , R , r, \omega,p,k} \ge\tau_{\ve , R , r,p}$ such that for every $\tau\ge \tau_{\ve , R , r, \omega,p,k}$ the  following holds:\\

The rescaled truncated \si-pencil $\tau \cdot \TP(p,\omega,R/\tau,r/\tau)\subset \tau \cdot H$ can be written* as a smooth section $\Gamma_\tau  \mm{ with } |\Gamma_\tau|_{C^k} < \ve$ of the  normal bundle of $C_p^\tau$ in $\R^{n+1}$ over $\TP(0,\omega,R,r)\subset C_p^\tau$.
\end{itemize}
\end{proposition}

\begin{remark}[Boundary adjustments]\label{ba}
 *Strictly speaking, $\Gamma_\tau$ slightly deviates from $\TP$ near $\p \TP$ . This deviation  disappears for $\tau \ra \infty$. Since we will only be interested in subsets situated within a bounded distance from the boundaries we are always free to adjust our definitions near the boundary according to our needs. We thus omit the precise convergence near the boundary. \qed
\end{remark}

\noindent{\bf Proof of Proposition~\ref{freez}} \, The \textbf{flat norm version} follows essentially from \ref{flat-norm-approx}. Indeed, assume that the required $\tau_{\omega , R , r,p}$ does not exist. Hence, there would be a sequence $\tau_i \ra \infty$, $i \ra \infty$ such that any $\tau_i \cdot (H \cap B_{R/\tau_i}(p)\setminus B_{r/\tau_i}(p))$ is not $\ve$-close in flat norm to any tangent cone. However, $\tau_i$ is subconvergent to some tangent cone $C$ by Proposition \ref{flat-norm-approx}. Hence there exists a sufficiently large $i_0 $ with $\tau_{i_0} \cdot (H \cap B_{R/\tau}(p) \setminus B_{r/\tau} (p))$ is $\ve$-close to $C \cap B_R(0) \setminus B_r (0 )$, contradicting that none of the  $\tau_i \cdot (H \cap B_{R/\tau_i}(p)\setminus B_{r/\tau_i}(p))$  is $\ve$-close to any tangent cone. \\

For the \textbf{$C^k$-version} 
assume again we could not find a $\tau_{\ve , R , r,\omega,p,k}$ as asserted. Then there would be a sequence $\tau_i \ra \infty$, $i \ra \infty$ such that ${\tau_i} \cdot H$ can never be written locally as a section $\Gamma_{\tau_i}$ with $|\Gamma_{\tau_i}|_{C^k} < \ve$ of the normal bundle over $\TP(0,\omega,R,r)\subset C_p$ for some tangent cone $C_p\in\mathcal{SC}_n$. However, 
by the compactness of $\mathcal{SC}_n$ from Proposition~\ref{cone-compactness} and the flat norm version just established, the sequence of truncated pencils $\TP(0,\omega,R,r) \subset C^{\tau_i}_p$ subconverges in flat norm to $\TP(0,\omega,R,r)\subset C$ for some tangent cone $C \in\mathcal{SC}_n$. By the $C^l$-norm convergence upgrade from Proposition~\ref{flat-norm-approx}, this implies subconvergence of $\tau_i\cdot H$ in $C^k$-norm, contradicting the assumption.\qed

%
%
\subsubsection{Non-Tangential Analysis and Extension Results}\label{nt00}
For many operators $L$, such as the conformal Laplacian or the Jacobi field operator, there is little control over positive solutions of $L\, f=0$ when we approach the singular set $\Sigma \subset H$ by classical means. The solutions develop poles, and the pole order can change with the local dimension of $\Sigma$ which is subject to jumping phenomena. By contrast, the quotient of any two positive solutions is remarkably well-behaved, at least if we approach $\Sigma$ \emph{non-tangentially} along \si-pencils $\P(p,\omega)$ pointing to any $p \in \Sigma$.\\

Towards this end we use the unfolding correspondence.

\begin{lemma}[\si-Pencils viewed relative $d_{\bp}$]\label{cy}
For any  $\omega >0$ and any geodesic ray $\gamma_p \subset (H \setminus \Sigma,d_{\bp})$, representing some $p \in \Sigma$, starting from a basepoint $p_0 \in H \setminus \Sigma$, and any sufficiently small $\eta >0$ there exists some $\zeta>0$ such that
\[\P(p,\omega) \cap B_{\eta}(p) \subset (H \setminus \Sigma, g_H) \mm{ is contained in the $\zeta$-distance tube } U_\zeta(\gamma_p) \subset  (H \setminus \Sigma,d_{\bp}) \mm{ of }\gamma_p.\]

Since $(H \setminus \Sigma, d_{\bp})$ and $(H \setminus \Sigma, d_{\bp^*})$ are quasi-isometric, this also gives an inclusion in the distance tube $U_{\zeta^*}(\gamma^*_p) \subset  (H \setminus \Sigma,d_{\bp^*})$, around a geodesic ray $\gamma^*_p \subset (H \setminus \Sigma,d_{\bp^*})$,  for some  $\zeta^* \ge  \zeta$.
\end{lemma}
\noindent\textbf{Proof} \, We  start from a given $\gamma_p$ in $(H \setminus \Sigma, d_{\bp})$. Recall that $\gamma_p$ is a c-\si-uniform curve relative to $(H \setminus \Sigma, g_H)$, for some $c > 0$. In particular, still relative $(H \setminus \Sigma, g_H)$, we have $l_{min}(\gamma_p(q)) \le c \cdot \delta_{\bp}(q)$ for any point $q \in \gamma_p$.
Since we are only interested in what happens close to $p$, we may assume that $l_{min}(\gamma_p(q))$ equals the length of the subarc to $p$, whence $d_H(q,p) \le c \cdot \delta_{\bp}(q)$. In particular, we may assume that $\gamma_p \subset \P(p,\omega)$ for any $\omega$ with $1/c > \omega >0$.\\

Now we employ the asymptotic \si-freezing \ref{freez} for $H$ around $p$: For given $\ve > 0$, $1/(2 \cdot c) > \omega>0$ and $R >1 >   r >0$ we have some $\tau_{\ve , R , r, \omega,p,5}> 1$ so that for any $\tau \ge \tau_{\ve , R , r, \omega,p,5}$, there is a tangent cone $C_p^\tau$ of $H$ at $p$ so that the subset $
\tau \cdot  B_{R/\tau}(p) \setminus B_{r/\tau} (p) \cap   \P(p,\omega)  \subset \tau \cdot H$
can be written as a smooth section $\Gamma_\tau  \mm{ with } |\Gamma_\tau|_{C^5} < \ve$  of the normal bundle of $B_R(0)\setminus B_r (0) \cap   \P(0,\omega)  \subset C_p^\tau.$\\

Therefore, choosing $\ve>0$ small enough, we observe that it is enough to understand the case where $\gamma_p \subset C_p^\tau$ with
\[
B_R(0)\setminus B_r (0) \cap \gamma_p \subset B_R(0)\setminus B_r (0) \cap   \P(0,\omega)
\]
for $1/(2 \cdot c)  > \omega >0$. Now we apply \cite[Proposition 2.7]{L1} and  \cite[Corollary 2.10]{L1} saying that any Euclidean minimal cone $C \in \mathcal{SC}_{n}$ is \si-uniform for a common  \si-uniformity constant $c_n\ge 1$ and $diam(S_C) \le d_n$ for some $d_n >0$ depending only on the dimension, where we set $S_C:=\p B_1(0) \cap C$. Then any two $v,w \in \{x \in S_C \setminus \Sigma_{S_C} \,|\, \delta_{\bp}(x) \ge \omega\} \subset S_C$ can be linked by an  \si-uniform curve
\[\gamma_{v,w} \subset \{x \in  C \setminus \Sigma_{C} \,|\, \delta_{\bp}(x) \ge \eta(\omega) \cdot \omega\}, \mm{ for some }\eta(\omega) \in (0,1),\mm{ and of length }\le l_n.\]
with both $\eta(\omega)$ and $l_n$ independent of $C$. Namely, using $\delta_{\bp}(x) \le L  \cdot dist(x,\Sigma)$  around the endpoints $v,w$ we have a radius $\rho(\omega)>0$  so that $\delta_{\bp}(z) \ge 2\omega/3$, for $z \in B_\rho$. Outside these balls the twisted double \si-cone estimate $l_{min}(\gamma_{v,w}(z)) \le c \cdot \delta_{\bp}(z)$ for any $z \in \gamma_{p,q}$ ensures the lower bound on $\delta_{\bp}$ along $\gamma_{v,w}$.  This entails the estimate
\[
d_{\bp}(y, B_R(0)\setminus B_r (0) \cap \gamma_p) \le l_n/(\eta(\omega) \cdot \omega) \mm{ for any } y \in S_C \cap   \P(0,\omega).
\]
The scaling property  $\bp_{\lambda \cdot H} \equiv \lambda^{-1} \cdot  \bp_{H}$, for any $\lambda >0$, then shows that for $\zeta:= 2l_n/(\eta(\omega) \cdot \omega)$ we have: $\P(p,\omega) \cap B_{\ve}(p)$  is contained in $U_\zeta(\gamma_p)$. \qed

Now we reach the counterpart of the classical Fatou theorem on $B_1(0) \subset \R^2$. We recall that  for \si-adapted operators any positive solution of $L\, f=0$ on $H \setminus \Sigma$ can be written in terms of the Martin integral $u_{\mu}(x) =\int_{\widehat{\Sigma}} k(x;y) \, d \mu(y)$ for some suitable finite Radon measure $\mu$.

\begin{theorem}[Relative Fatou Theorem on $\mathbf{H \setminus \Sigma}$]\label{ftt}
Let $H \in {\cal{G}}$ and $L$ be an \si-adapted operator on $H \setminus \Sigma$. Further, let $\mu$ and $\nu$ be two finite Radon measures on $\Sigma$ with associated solutions $u_{\mu}$ and $u_{\nu}$ of $L \, f=0$ as in \eqref{muu}. Then for $\nu$-almost any $p \in \Sigma$ and any fixed $\rho >0$, we have
\[
u_{\mu}/u_{\nu}(x) \ra d \mu/d \nu(p)\,\mm{ as } x \ra p, \mm{ with } x \in \P(p,\omega).
\]
(Here, $d \mu/d \nu$ denotes the Radon-Nikodym derivative of $\mu$ with respect to $\nu$.)
\end{theorem}

To explain this statement we recall that for any two Radon measures $\mu$, $\nu$ on $\Sigma$ we have with respect to $\nu$ a (uniquely determined) Lebesgue decomposition $\mu=\mu_1 + \mu_2$ into a $\nu$-absolute continuous measure $\mu_1$ and a $\nu$-singular measure $\mu_2$. Concretely, there exists a $\nu$-integrable function $f$, the so-called \textbf{Radon-Nikodym derivative} $d \mu_1/d \nu$ such that $\mu_1(E)= \int_E f \, d\nu$ for any measurable $E \subset \Sigma$. Note that $f$ is uniquely determined up to $\nu$-negligible sets. On the other hand, there exists a $\nu$-negligible set $F \subset \Sigma$ such that $\mu_2(\Sigma \setminus F)=0$. We can thus rephrase Theorem~\ref{ftt} as follows. For any given $\rho >0$ we have for $\nu$-almost any $p \in \Sigma$ that
\[
u_{\mu}/u_{\nu}(x) \ra f(p) = d \mu_1/d \nu(p)\,\mm{ as } x\in \P(p,\omega) \ra p.
\]

\noindent\textbf{Proof} \,
By \cite[Theorem 4]{An1} (or \cite[Theorem 6.5 p.\ 100 and definitions on p.\ 99]{An2}) which, in turn, is an adaptation of the work by Gowrisankaran~\cite{Go}, the unfolding correspondence in the proof of \ref{mbhs} yields the following result for $(H \setminus \Sigma, d_{\bp^*})$ equipped with
$\delta_{\bp^*}^2 \cdot L$ : For any two finite Radon measures $\mu$, $\nu$ on $\p^0_M((H \setminus \Sigma, d_{\bp^*}),\delta_{\bp^*}^2 \cdot L)$ we get for almost any $p \in \p^0_M((H \setminus\Sigma, d_{\bp^*}),\delta_{\bp^*}^2 \cdot L)$ and any distance tube $U_r(\gamma_p)$, $r >0$, around a geodesic ray $\gamma_p \subset ((H \setminus \Sigma,d_{\bp^*})$ representing $p$, that the solutions $u_{\mu}$, $u_{\nu}$ associated with $\mu$, $\nu$ satisfy
\[
u_{\mu}/u_{\nu}(x) \ra d \mu/d \nu(p)\,\mm{ as } x \ra p, \mm{ with } x \in U_r(\gamma_p).
\]
However, Lemma~\ref{cy} shows that for any $\rho >0$ and sufficiently small $\ve >0$ there exists (a large) $r >0$ so that $\P(p,\omega) \cap B_{\ve}(p)$ is contained in $U_r(\gamma_p)$. The result follows.  \qed

\begin{remark}
For measures with positive density along $\Sigma$ there is actually no general control over the tangential behavior of the quotients $u_{\mu}/u_{\nu}$. In fact, examples due to Littlewood which were refined by Aikawa, cf.\ \cite{Ai4}, show that bounded harmonic functions on the unit disc in $\R^2$ exists so that the limit along paths approaching $S^1$ tangentially does \emph{not} exist in any $z \in S^1$. It is conceivable that similar examples also exist for \si-adapted operators on $H \setminus \Sigma$.\qed
\end{remark}

Next we turn to a situation unaddressed by Fatou's theorem. Namely, we consider quotients of solutions $u_{\mu}$, $u_{\nu} >0$ on $H \setminus \Sigma$ whose associated Radon measures $\mu$, $\nu$ vanish on a common open subset $A \subset \Sigma$ with $\mu(A)=\nu(A)=0$. Still, by virtue of the BHP, the quotient $u_{\mu}/u_{\nu}$ remains well-controlled even in tangential directions. To see this we first prove the

\begin{lemma}\label{cea}
Let $H \in {\cal{G}}$ and $L$ be an \si-adapted operator on $H \setminus \Sigma$. For any open subset $A \subset \widehat{\Sigma}$, $ \mu(A)=0$ implies that  $u_{\mu}$ $L$-vanishes along $A$.
\end{lemma}
\noindent\textbf{Proof} \,
If $\mu(A)=0$, then the Martin integral and the inequality \eqref{fhepq1} for minimal functions coming from BHP for Green's Functions, show the following. For each $z \in A$ there exists a constant $c_z>0$ such that for the Green's functions $G(\cdot,p)$ relative to some basepoint $p \in H \setminus \Sigma$,
\begin{equation}\label{lva}
G(\cdot,p) \ge c_z \cdot u_{\mu}\mm{ near }z.
\end{equation}
But $G(\cdot,p)$ $L$-vanishes along $\Sigma$ and thus $u_{\mu}$ $L$-vanishes along $A$.\qed

The latter result shows that \emph{L}-vanishing of solutions along $A$ is naturally related to the vanishing of the associated Radon measure they define. For any pair of such solutions we can now formulate a remarkable counterpart of the Fatou theorem.

\begin{theorem}[Continuous Extensions to $\mathbf{\Sigma}$]\label{cee}
Let $H \in {\cal{G}}$ and $L$ be an \si-adapted operator on $H \setminus \Sigma$. For any two solutions $u$, $v>0$ of $L\, f=0$,  both $L$-vanishing along some common open subset $A \subset \widehat{\Sigma}$, the quotient $u/v$ on $H \setminus \Sigma$ admits a continuous extension to $H \setminus \Sigma \cup A \subset \widehat{H}$.
\end{theorem}
\noindent\textbf{Proof} \,
We modify a classical method due to Moser in \cite[Chapter 5]{Mo}, cf.\ \cite[Theorem 7.9]{JK} and \cite[Theorem 2]{Ai2} to the case of boundary problems. Originally, this theory was employed to derive relative estimates for the oscillation of harmonic functions on concentric Euclidean balls.\\

Instead of concentric balls we shall use a canonical $\Phi_\delta$-chain $\mathcal{N}^\delta_k(\gamma)$ along a geodesic ray representing $z \in \widehat{\Sigma}$ cf.\ref{cpc}. As a result we get Harnack estimates on the $\mathcal{N}^\delta_k$ from the BHP for \si-adapted operators in  ~\ref{mbhsq}. We set
\begin{equation}\label{mi}
\sup(k):= \sup_{\mathcal{N}^\delta_k} u/v  \mm{ and } \inf(k):= \inf_{\mathcal{N}^\delta_k} u/v.
\end{equation}
The \emph{oscillation} of $u/v$ on $\mathcal{N}^\delta_k$ is the quantity $osc(k)= \sup(k)-\inf(k)$.\\

We first note from Proposition~\ref{mbhsq} that $\sup(k_0) < \infty$ for any $k_0$ sufficiently large. Hence the poles of $u$ and $v$ do not belong to $\mathcal{N}^\delta_{k_0-2}$. For $k \ge k_0$, we consider the two solutions $\sup(k) \cdot v -u$ and  $u - \inf(k) \cdot v$. They are positive by Hopf's maximum principle and $L$-vanish along $A$. Therefore the BHP in  ~\ref{mbhsq} applies to this pair of functions and shows that for $\mathcal{N}^\delta_{k+1}$,
\begin{align*}
\sup_{\mathcal{N}^\delta_{k+1}} (\sup(k) \cdot v -u)/v  &\le C^* \cdot  \inf_{\mathcal{N}^\delta_{k+1}} (\sup(k) \cdot v -u)/v\\
\sup_{\mathcal{N}^\delta_{k+1}} (u - \inf(k) \cdot v)/v  &\le C^* \cdot  \inf_{\mathcal{N}^\delta_{k+1}} (u - \inf(k) \cdot v)/v.
\end{align*}
From these inequalities we get
\[
\sup(k) - \inf(k+1) \le C^*  \cdot\big(\sup(k) - \sup(k+1)\big) \,\mm{ and } \,\sup(k+1) - \inf(k) \le C^*  \cdot \big(\inf(k+1) - \inf(k)\big).
\]
We add suitable multiples of these inequalities so that
\begin{align*}
osc(k+1)&= \sup(k+1) - \inf(k+1) \le(C^*-1)/(C^*+1) \cdot (\sup(k) - \inf(k))\\
&=(C^*-1)/(C^*+1) \cdot osc(k).
\end{align*}
Hence, setting $a:=(C^*-1)/(C^*+1) <1$,
\begin{equation}\label{oscc}
osc(k) \le a^{k-k_0} \cdot osc(k_0) \ra 0 \mm{ as } k \ra \infty.
\end{equation}
Consequently, $u/v$ remains bounded near $z$ as a result from the BHP. By \eqref{oscc}, the quotient admits a continuous extension to $z$.\qed

For certain types of operators we also yield solvability of the Dirichlet problem:

\begin{theorem}[Dirichlet Problem for \si-Adapted Operators]\label{diri}
Let $H \in {\cal{G}}$ and $L$ be an \si-adapted operator on $H \setminus \Sigma$ such that
\begin{itemize}
  \item constant functions solve $L \, f=0$;
  \item for a given $p \in H \setminus \Sigma$, the Green's function $G(x,p)\ra 0$ as $x \ra \widehat{\Sigma}$.
\end{itemize}
Then, for any continuous function $f$ on $\widehat{\Sigma}$, there exists a uniquely determined continuous function $F$ on $H$ such that
\[
L \, F =0 \mm{ with } F|_{\widehat{\Sigma}} \equiv f.
\]
\end{theorem}
\noindent\textbf{Proof} \,
The assumptions and the Martin theory for \si-adapted operators allow us to imitate the standard arguments which are, for instance, well-known for harmonic functions on Euclidean domains.\\

We first let $\mu^{\cs}$ be the Radon measure on $\widehat{\Sigma}$ which is associated with the constant function $1$. Then we decompose $f$ into $f=f^+-f^-$, for $f^+:=\max\{f,0\}$, $f^-:=-\min\{f,0\}$ and set
\[
F^\pm(x) =  \int_{\widehat{\Sigma}} f^\pm(y) \cdot k(x;y) \,d \mu^{\cs}(y).
\]
It follows that $F^\pm \ge 0$ solves $L \, f=0$. We assert that $F^\pm$ extends continuously to $H$
and coincides with $f^\pm$ on $\widehat{\Sigma}$. This gives the claim for $F:=F^+-F^-$.\\

For any $z \in \widehat{\Sigma}$ we have to show that $F^\pm(x)\ra f^\pm(z)$ as $x \ra z$, $x \in H \setminus \Sigma$. Towards this end we notice from \eqref{lva} that for any fixed neighborhood $U \subset \widehat{\Sigma}$ of $z$,
\[
G(\cdot,p) \ge c_z \cdot \int_{\widehat{\Sigma} \setminus U} f^\pm(y) \cdot k(x;y) \,d \mu^{\cs}(y)
\]
near $z$. Thus, for any such $U$, we have $\int_{\widehat{\Sigma} \setminus U} f^\pm(y) \cdot k(x;y) \,d \mu^{\cs}(y) \ra 0$ as $x \ra z$.\\

On the other hand, $f^\pm$ is continuous on $\widehat{\Sigma}$, so $|f^\pm(y) - f^\pm(z)| < \ve$ for $y \in U(\ve)$, where $U(\ve) \subset \widehat{\Sigma}$ is a sufficiently small neighborhood of $z$. Moreover, since $ \int_{\widehat{\Sigma}}  k(x;y) \,d \mu^{\cs}(y)=1$, the definition of $\mu^{\cs}$ and $\int_{\widehat{\Sigma} \setminus U(\ve)} k(x;y) \,d \mu^{\cs}(y) \ra 0$ as $x \ra z$ imply that for any $x \in H \setminus \Sigma$ close enough to $z$,
\[
\left|\int_{U(\ve)} f^\pm(y) \cdot k(x;y) \,d \mu^{\cs}(y) - \int_{U(\ve)} f^\pm(z) \cdot k(x;y) \,d \mu^{\cs}(y)\right| \le 2\ve \mm{ for } y \in U(\ve).
\]
We conclude that $F^\pm(x)  \ra f^\pm(z)$ as $x \ra z$. Finally, the maximum principle shows that $F \equiv 0$ is the only solution with $F|_{\widehat{\Sigma}} \equiv 0$. This yields the asserted unique solvability of the Dirichlet problem.\qed

%
%
%
\setcounter{section}{5}
\renewcommand{\thesubsection}{\thesection}
\subsection{Criticality and Principal Eigenvalues}
In this chapter we mainly focus on eigenvalue problems. It is therefore convenient to restrict to the case of symmetric operators where we can exploit variational arguments.\\

Actually, the eigenvalue theory for many $\bp$-adapted operators which do \emph{not} satisfy the weak coercivity condition can be discussed by considering the \si-adapted operators $L_\lambda:=L-\lambda \cdot \bp^2 \cdot Id$ for suitable $\lambda \in \R$. This is an important tool for extending the range of problems our Martin theory can be applied to.\\

We first introduce \si-Sobolev spaces which we shall use on several occasions in this and the follow-up paper~\cite{L2}. For instance, they are crucial in ~\cite[Theorem 3 and 5]{L2}. We believe they prove useful for many other purposes as well.

%
%
\subsubsection{\si-Sobolev Spaces}\label{ssa}
In presence of an \si-structure $\bp$ we can define the following natural Hilbert space structures on subspaces of $M(H \setminus \Sigma)$, the space of measurable functions on $H \setminus \Sigma$. Like for most results employing \si-structures the minimal requirement for the following concepts is that $H$ is \textbf{not} totally geodesic and hence $\bp>0$. Recall that totally geodesic $H\in\cal{G}$ are automatically smooth (for  ${\cal{H}}^{\R}_n$ they are just hyperplanes) since $|A|$ diverges towards singularities.

\begin{definition}\label{a}
For any $f$, $g \in M(H \setminus \Sigma)$, $H\in\cal{G}$  not totally geodesic, we define
 \begin{itemize}
  \item the $L^2_{\bp}$-norm:\, $|f|_{L^2_{\bp}(H \setminus \Sigma)}:=(\int_{H \setminus \Sigma} \bp^2 \cdot f^2 \, dV)^{1/2}$.
  \item  the $L^2_{\bp}$-scalar product:\, $\langle f,g \rangle_{L^2_{\bp}(H \setminus \Sigma)}:= \int_{H \setminus \Sigma} \bp^2 \cdot f \cdot g \, dV$.
  \item the space $L^2_{\bp}(H \setminus \Sigma)$ of $L^2_{\bp}$-integrable functions
  \begin{equation}\label{laa}
  L^2_{\bp}(H \setminus \Sigma):=\{f \in M(H \setminus \Sigma) \,|\ |f|_{L^2_{\bp}(H \setminus \Sigma)} < \infty\}.
  \end{equation}
\end{itemize}
If there is no risk of confusion we shall write $L^2_{\bp}$ for short. Also, we define Sobolev norms and scalar products for functions $f$, $g \in M(H \setminus \Sigma)\cap L^1_{loc}(H \setminus \Sigma)$ as follows:
\begin{itemize}
  \item  the $H^{1,2}_{\bp}$-norm:\,  $|f|_{H^{1,2}_{\bp}(H \setminus \Sigma)}:=(\int_{H \setminus \Sigma}|\nabla f|^2 + \bp^2 \cdot f^2 \, dV)^{1/2}$  where $\nabla f$ is the distributional derivative.
  \item the $H^{1,2}_{\bp}$-scalar product:\, $\langle f,g \rangle_{H^{1,2}_{\bp}(H \setminus \Sigma)}:= \int_{H \setminus \Sigma} \langle \nabla f, \nabla g \rangle + \bp^2 \cdot f \cdot g \, dV $.
  \item the $H^{1,2}_{\bp}$-integrable, measurable functions in terms of distributional derivatives:
  \begin{equation}\label{lab}
  H^{1,2}_{\bp}(H \setminus \Sigma):=\{f \in M(H \setminus \Sigma)\,|\ |f|_{H^{1,2}_{\bp}(H \setminus \Sigma)} < \infty\}.
  \end{equation}
  Equivalently, this is the $H^{1,2}_{\bp}(H \setminus \Sigma)$-completion of the space of smooth functions $C^\infty(H \setminus \Sigma)\cap H^{1,2}_{\bp}(H \setminus \Sigma)$. Again, we shall simply write $H^{1,2}_{\bp}$.
\end{itemize}
\end{definition}

By standard arguments, $L^2_{\bp}(H \setminus \Sigma)$ and $H^{1,2}_{\bp}(H \setminus \Sigma)$ are Hilbert spaces.\\

\begin{remark}
Note that $\bp_{\lambda\cdot H}=\lambda^{-1}\cdot\bp_H$ has the same scaling behaviour as the dual metric $\langle\cdot,\!\cdot\rangle$ on $1$-forms. Hence both integrands in the $H^{1,2}_{\bp}$-norm transform in the same way under scalings of the underlying space.
\end{remark}

This has the following valuable consequence.  Denote the space of smooth functions with compact support on $H \setminus \Sigma$ by $C^\infty_0(H \setminus \Sigma)$.

\begin{theorem}[Compactly Supported Approximations]\label{ba}
Let $H\in\cal{G}$ be a non totally geodesic hypersurface. Then we have
\begin{equation}\label{csa}
H^{1,2}_{\bp}(H \setminus \Sigma) \equiv H^{1,2}_{\bp,0}(H \setminus \Sigma):= H^{1,2}_{\bp}\mm{-completion of }C^\infty_0(H \setminus \Sigma).
\end{equation}
\end{theorem}

Before we start with the proof we define for present and later use some suitable cut-off functions concentrated near $\Sigma$. We fix some $\psi \in C^\infty(\R,[0,1])$ with $\psi \equiv 1$ on $\R^{\le 0}$ and $\psi \equiv 0$ on $\R^{\ge 1}$. Then we set
\begin{equation}\label{cut}
\psi[\eta](x):= \psi(\eta^{-1} \cdot \delta_{\bp}(x)-1)\mm{ for } x \in H \setminus \Sigma \mm{ and } \eta \in (0,1).
\end{equation}
Recall from (S3) in Definition~\ref{def1} that $\delta_{\bp}$ is Lipschitz continuous:
\[
|\delta_{\bp}(p)- \delta_{\bp}(q)|   \le L_{\bp} \cdot d(p,q) \mm{ for } p,q \in  H \setminus \Sigma \mm{ and some constant } L_{\bp}>0.
\]
Hence, upon applying the chain rule, there exists some constant $c(\psi) >0$ depending only on the chosen cut-off function $\psi$ such that
\begin{equation}\label{nabb}
|\nabla \psi[\eta](x)| \quad
\begin{cases}
\le c(\psi) \cdot \bp(x), &\text{for $\delta_{\bp}(x) \in (\eta, 2 \cdot \eta)$;}\\
= 0, &\text{otherwise.}
\end{cases}
\end{equation}
(Here, we either interpret the derivatives of $\delta_{\bp}$ distributionally or use Rademacher's theorem. Alternatively, one can apply the Whitney smoothing \cite[Proposition B.3]{L1} and use $\delta^*_{\bp}$.) \\

\noindent\textbf{Proof} \,
Let $f \in H^{1,2}_{\bp}(H \setminus \Sigma)$.\\

Assume first that $H$ is compact. When $H$ is regular, the claim is trivial, so we may assume that $\Sigma_H \n$. For $1>\eta>0$ we consider the functions $(1-\psi[\eta]) \cdot f$ which are compactly supported in $H \setminus \Sigma$. Further, we set $\chi_{\eta}(x)=1$ for $\delta_{\bp}(x) \in (\eta, 2\eta)$ and $\chi_{\eta}(x)=0$ otherwise. We claim that
\[
(1-\psi[\eta]) \cdot f \ra f \mm{ as } \eta \ra 0 \mm{ in } H^{1,2}_{\bp}\mm{-norm}.
\]
Indeed, from \eqref{nabb} and $\chi_{\eta}$, $\psi[\eta] \le  \psi[2\eta]$ we get
\begin{align*}
|(1- (1-\psi[\eta])) \cdot f|^2_{H^{1,2}_{\bp}}&= \int_{H \setminus \Sigma}|\nabla (\psi[\eta] \cdot f)|^2 + \bp^2 \cdot (\psi[\eta] \cdot f)^2 \, dV\\
&\le 2\int_{H \setminus \Sigma}\psi[\eta]^2 \cdot | \nabla f|^2 + |\nabla \psi[\eta]|^2\cdot f^2 + \bp^2 \cdot \psi[\eta]^2 \cdot f^2 \, dV\\
&\le 2\int_{H \setminus \Sigma}\psi[\eta]^2 \cdot | \nabla f|^2 + (c(\psi)^2 \cdot  \chi^2_{\eta}(x) + 1) \cdot \bp^2 \cdot \psi[\eta]^2 \cdot f^2 \, dV\\
&\le 2(c(\psi)^2  + 1) \cdot \int_{H \setminus \Sigma}\psi[2\eta]^2 \cdot | \nabla f|^2 + \bp^2 \cdot \psi[2\eta]^2 \cdot f^2 \, dV  \\
&\le 2(c(\psi)^2  + 1) \cdot \int_{\{ x \in H \setminus \Sigma \, |\, \delta_{\bp}(x) \le 2\eta\}}  | \nabla f|^2 + \bp^2 \cdot  f^2 \, dV.
\end{align*}
But the latter term tends to $0$ as $\eta\to0$, for $|f|_{H^{1,2}_{\bp}(H \setminus \Sigma)} < \infty$ and $Vol(\{ x \in H \setminus \Sigma \, |\, \delta_{\bp}(x) \le 2 \cdot \eta\}) \ra 0$ as $\eta \ra 0$ by compactness of $H$.\\

For non-compact $H \in \cal{G}$, i.e., $H$ is a minimal boundary in Euclidean space $H \in {\cal{H}}^{\R}_n$, we introduce a cut-off towards infinity by setting $\psi_R(x):= \psi(R^{-1} \cdot |x|-1)$ for $x \in H \setminus \Sigma$ and $R \ge 1$. We show that for any $f \in H^{1,2}_{\bp}(H \setminus \Sigma)$,
\[
\psi_R \cdot f \ra f \mm{ as } R \ra \infty \mm{ in } H^{1,2}_{\bp}\mm{-norm}.
\]
Toward this end we note that
\begin{equation}\label{nabba}
|\nabla \psi_R(x)| \quad
\begin{cases}
\le c(\psi)^*/R, &\text{for $|x| \in (R, 2 \cdot R)$;}\\
= 0, &\text{otherwise}
\end{cases}
\end{equation}
for some $c(\psi)^* >0$ depending only on the chosen $\psi$. Further, we let $\chi_R^*(x)=1$, $x \in B_{2R}(0) \setminus B_R(0)$, and $\chi_R^*(x)=0$ if otherwise.  The Lipschitz condition $  |\delta_{\bp_H}(p)- \delta_{\bp_H}(q)|   \le L_{\bp} \cdot d_H(p,q) $ applied to some basepoint $p=p_0 \in H$ and a variable point $q=x$ shows that $\delta_{\bp_H}(x)  \le L_{\bp} \cdot d_H(x,p_0) + \delta_{\bp_H}(p_0)$. This means  that $\bp(x) \ge a_H \cdot d_H(x,p_0) ^{-1}$ on $H \setminus (\Sigma_H \cup B_1(p_0))$ for some $a_H>0$. Then we have for $R >1$:

\begin{align*}
|(1-\psi_R) \cdot f|^2_{H^{1,2}_{\bp}}&= \int_{H \setminus \Sigma}|\nabla ((1-\psi_R) \cdot f)|^2 + \bp^2 \cdot ((1-\psi_R) \cdot f)^2 \, dV\\
&\le 2 \cdot \int_{H \setminus \Sigma}(1-\psi_R)^2 \cdot | \nabla f|^2 + |\nabla \psi_R|^2\cdot f^2 + \bp^2 \cdot (1-\psi_R)^2 \cdot f^2 \, dV\\
&\le 2 \cdot \int_{H \setminus \Sigma}(1-\psi_R)^2 \cdot | \nabla f|^2 + \bp^2 \cdot (1-\psi_R)^2 \cdot f^2 \, dV \\
&\quad  +  2 \cdot \int_{B_{2 \cdot R}(p_0)\setminus B_R(p_0)}  (c(\psi)^*/R)^2 \cdot f^2 \, dV\\
&\le  \left(2 +(c(\psi)^*/a_H)^2\right) \cdot \int_{(H \setminus \Sigma) \setminus B_R(p_0)}  | \nabla f|^2 + \bp^2 \cdot  f^2 \, dV,
\end{align*}
and the latter terms tends to $0$ as $R \ra \infty$. Combining this with the argument for the compact case applied to $\overline{B_R} \cap H$, we get the claim. \qed

%
%
\subsubsection{Hardy Inequalities}
Next we consider $\bp$-adapted operators on $H \setminus \Sigma$  which are symmetric for the standard $L^2$-inner product. For these the $\bp$-weak coercivity can be expressed in terms of a variational integral for eigenvalues. This is a well-known coercivity property of bilinear forms, cf.\ \cite[Proposition 1 and Appendix]{An3} for an explicit treatment of the Laplacian on Euclidean domains. For completeness, we explain how this carries over to our case.

\begin{theorem}[Hardy Inequalities]\label{hi}
Let $L$ be an $\bp$-adapted symmetric operator on $H \setminus \Sigma$. Then $L$ is $\bp$-weakly coercive if and only if there exists a positive constant $\tau = \tau(L,\bp,H)>0$ such that the \textbf{Hardy inequality}
\begin{equation}\label{hadi}
\int_{H \setminus \Sigma}  f  \cdot  L f  \,  dV \, \ge \, \tau \cdot \int_{H \setminus \Sigma}\bp^2\cdot f^2 dV
\end{equation}
holds for all $f \in C^\infty_0(H \setminus \Sigma)$.
\end{theorem}

\begin{remark}\label{h}
Consider the symmetric quadratic form associated with $L$ defined by $b(u,v):= \int_{H \setminus \Sigma} u \cdot L v$, $u$, $v\in C^\infty_0(H \setminus \Sigma)$. Via partial integration we can rewrite $b$ as a Dirichlet type integral in terms of $u$, $v$ and their first derivatives. Since $L$ is $\bp$-adapted, $b(u,v)$ admits a canonical extension to $H^{1,2}_{\bp,0}(H \setminus \Sigma)=H^{1,2}_{\bp}(H \setminus \Sigma)$ satisfying \eqref{hadi} for the same constant $\tau>0$. Furthermore, a supersolution $f \in H^{1,2}_{\bp}(H \setminus \Sigma)$ satisfies $L f > \ve \cdot \bp^2 \cdot f $ in the weak sense, that is, $b(f,v) \ge 0$ for each nonnegative test function $v \in C^\infty_0(H \setminus \Sigma)$. This matches the condition from Definition~\ref{sao0} used in Ancona's theory.\qed
\end{remark}

\noindent\textbf{Proof} \,
When \eqref{hadi} is satisfied we consider the continuous quadratic form
\[
a(u,v):= \int_{H \setminus \Sigma} (u \cdot L v  - \ve \cdot  \bp^2 \cdot u \cdot v) \, dV \mm{ on } C^\infty_0(H \setminus \Sigma).
\]
From Lemma \ref{ba} we can extend $a(\cdot,\cdot)$ to $H^{1,2}_{\bp}(H \setminus \Sigma)$. For $\ve/ \tau \in (0,1)$ we then have
\[
a(f,f) \ge \big(1- \frac{\ve}{\tau}\big) \int_{H \setminus \Sigma} f\cdot L f \, dV \ge \big(1- \frac{\ve}{\tau}\big)/\big(1+ \frac{1}{\tau}\big) \cdot |f|_{H^{1,2}_{\bp}(H \setminus \Sigma)}.
\]
Thus, $a(\cdot,\!\cdot)$ is again a positive definite scalar product on $H^{1,2}_{\bp}(H \setminus \Sigma)$. For $h \in H^{1,2}_{\bp}(H \setminus \Sigma)$ we then consider the continuous functional $I_h: v \mapsto \int_{H \setminus \Sigma} h \cdot v \, dV$ on $H^{1,2}_{\bp}(H \setminus \Sigma)$. By Riesz' representation theorem there exists a unique $f=f(h) \in H^{1,2}_{\bp}(H \setminus \Sigma)$ with $a(f,v)= I_h(v) = \int_{H \setminus \Sigma} h \cdot v \, dV$ for all $v \in H^{1,2}_{\bp}(H \setminus \Sigma)$. From the definitions it follows that $L f - \ve \cdot \bp^2 \cdot f=h$. Now if we choose some smooth $h >0$ and $v=\min\{0,f\}$, then
\[
0 \le a(\min\{0,f\},\min\{0,f\})= a(f,\min\{0,f\})= \int_{H \setminus \Sigma} h \cdot \min\{0,f\}\, dV.
\]
Hence $f \ge 0$ and $L f > \ve \cdot \bp^2 \cdot f $. After some perturbation, small in $C^2$-norm, we get $f > 0$ which still satisfies $L f > \ve \cdot \bp^2 \cdot f $.\\

Now, if such an $f>0$ exists, we consider the first Dirichlet eigenvalue $\lambda_D$ of $L-\ve \cdot \bp^2 \cdot Id$ on a smoothly bounded domain $D$ with $\overline{D} \subset H \setminus \Sigma$ and its positive eigenfunction $w$, that is, $(L-\ve \cdot \bp^2 \cdot Id) w = \lambda_D \cdot w$. Then
\[
0=\int_D f \cdot (L -\ve \cdot \bp^2 \cdot Id- \lambda_D) \, w \, dV =\int_D w \cdot (L -\ve \cdot \bp^2\cdot Id - \lambda_D) \, f  \,dV  \ge -\lambda_D \int_D w \cdot   f  \,dV .
\]
Thus $\lambda_D \ge 0$ and $\int_{H \setminus \Sigma}  v \cdot  (L-\ve \cdot \bp^2 \cdot Id) v \ge 0$ for any $v  \in C_0^\infty(D)$ by standard extremal properties of the first eigenvalue (Rayleigh quotient). Since $D$ can be chosen arbitrarily large within $H \setminus \Sigma$, we get the Hardy inequality for $\tau = \ve$. \qed

\noindent\textbf{Shifted \si-Adaptedness} \,
By the previous lemma there exists for any symmetric \si-adapted operator $L$ a \emph{largest} $\tau >0$ such that the Hardy inequality \eqref{hadi} holds. This characterization of \si-adaptedness suggests to consider also the case where $\tau$ is negative but finite: $0>\tau > -\infty$.\\

The largest value for $\tau$ in \eqref{hadi} can be regarded as an eigenvalue of the operator $\delta_{\bp}^2 \cdot L$. Theorem \ref{scal2} will show that if $H$ is \emph{singular}, then for any value smaller than this largest value, positive eigenfunctions of $\delta_{\bp}^2 \cdot L$ do exist. This is an important difference to the case of smooth compact manifolds. Here, the first eigenvalue $\lambda_1$ of such an elliptic operator is the unique eigenvalue with an, up to constant multiples, unique positive eigenfunction, cf.\ \cite[Chapter VI]{C}. This will be a key feature for applications of singular area minimizers in scalar curvature geometry, cf.\ \cite[Theorem 2]{L2}.\\

Finally, we introduce the following terminology.

\begin{definition}\label{saop}
Let $L$ be an $\bp$-adapted symmetric operator on $H \setminus \Sigma$. If there is a \textbf{finite} $\tau \in \R$ such that the Hardy inequality \eqref{hadi} holds, then we call $L$ a \textbf{shifted \si-adapted} operator. We denote the largest such $\tau$ by $\lambda^{\bp}_{L,H}$ and call it the \textbf{(generalized) principal eigenvalue} of $\delta_{\bp}^2 \cdot L$.
\end{definition}

In particular, $L$ is \si-adapted if and only if $\lambda^{\bp}_{L,H} >0$.
%
%
\subsubsection{Basic Spectral Theory}\label{bgrb}
Let $H$ be either compact with $\Sigma \n$, or $H \subset \R^{n+1}$ be not totally geodesic. If $L$ is  shifted \si-adapted on $H\setminus\Sigma$, we get the following neat description of the spectral theory of $\delta_{\bp}^2 \cdot L$.

\begin{theorem}[Criticality and Principal Eigenvalues]\label{scal2}
Let $H \in \cal{G}$ and $\Sigma_H \n$, and $L$ be a shifted \si-adapted operator on $H \setminus \Sigma$. We set
\[
L_\lambda:= L - \lambda \cdot \bp^2 \cdot Id, \mm{ for }\lambda \in\R.
\]
Then we have the following trichotomy for the spectral theory of $\delta_{\bp}^2 \cdot L$.
\begin{itemize}
    \item \emph{\textbf{Subcritical}} when $\lambda < \lambda^{\bp}_{L,H}$: $L_\lambda$ is \si-adapted.
    \item \emph{\textbf{Critical}} when $\lambda = \lambda^{\bp}_{L,H}$:  There exists, up to multiples, a unique positive solution $\phi$ of $L _{\lambda^{\bp}_{L,H}} \ f = 0$, the so-called \textbf{ground state}\footnote{If $L$ is a Schr\"odinger operator, then $\phi$ \emph{is} the ground state well-known from the study of quantum mechanical systems, cf.\ \cite[Chapter 3.3]{GJ}.} of $L_{\lambda^{\bp}_{L,H}}$. The function $\phi$ can be described as the limit of first Dirichlet eigenfunctions for the operator $\delta_{\bp}^2 \cdot L$ on a sequence of smoothly bounded domains $\overline{D}_m \subset D_{m+1} \subset H \setminus \Sigma$, $m \ge 0$, with $\bigcup_m D_m = H \setminus \Sigma$.
    \item \emph{\textbf{Supercritical}} when $\lambda > \lambda^{\bp}_{L,H}$: $L_\lambda \, u = 0$ has no positive solution.
\end{itemize}
\end{theorem}

\begin{remark}\label{regu}
In the critical case it is easy to see that $L_\lambda$ does not admit a positive Green's function. This can be shown the same way we prove uniqueness of the ground state in the argument below. Following the conventions in \cite[Chapter 4]{P} we therefore say the operator $L$ is \begin{itemize}
  \item \textbf{subcritical}, if $L$ admits a positive Green's function.
  \item \textbf{critical} if it does \emph{not} admit a Green's function but a positive solution of $L \, f=0$.
  \item \textbf{supercritical} when the latter equation does \emph{not} admit any positive solutions.\qed
\end{itemize}

\end{remark}

\textbf{Proof of Theorem \ref{scal2}} \,
The main case is the critical case, where $\lambda = \lambda^{\bp}_{L,H}$. The subcritical case is rather obvious. The supercritical case follows from the discussion of the critical case.\\

\textbf{Subcritical Case} ($\lambda < \lambda^{\bp}_{L,H}$): \,
The variational relation $\int_H  f  \cdot L  f  \,  dV \, \ge \, \lambda^{\bp}_{L,H} \cdot \int_H \bp^2\cdot f^2 dV$ shows that
\[
\int_H  f \cdot L_{\lambda} f  \,  dV \, \ge\, (\lambda^{\bp}_{L,H} - \lambda) \cdot \int_H \bp^2\cdot f^2 dV.
\]
We conclude  that $L_{\lambda}$ is \si-adapted.\\

\textbf{Critical Case} ($\lambda = \lambda^{\bp}_{L,H}$) \,
We choose a sequence of smoothly bounded domains $\overline{D}_m \subset D_{m+1} \subset H
\setminus \Sigma$, $m \ge 0$, with $\bigcup_m D_m = H \setminus \Sigma$.\\

We consider the uniquely determined first eigenvalue $\lambda_m$ with eigenfunction $\phi_m$ on $\overline{D}_m$ for the Dirichlet problem of the operator $\delta_{\bp}^2 \cdot L$, that is,
\[
L  \, \phi_m = \lambda_m \cdot \bp^2 \cdot \phi_m \mm{ and } \phi_m \equiv 0\mm{ on } \p D_m.
\]
We may assume that $\phi_m > 0$ on $D_m$ and $\phi_m(p_0) = 1$ in some base point $p_0 \in \bigcap_m D_m$. Further, we formally extend $\phi_m$ to $H \setminus \Sigma$ by setting $\phi_m\equiv0$ outside $D_m$. The variational characterization of $\lambda_m$ via the Rayleigh quotient
\[
\lambda_m = \inf \Big\{\int_H  f  \cdot  L f  \,  dV / \int_H \bp^2\cdot f^2 dV \,\Big|\, f \mm{ smooth,
supported in } D_m \subset H \setminus \Sigma\Big\}
\]
shows that $\lambda_m > \lambda_{m+1} > \lambda^{\bp}_{L,H} >0$. Moreover, $\lim_{m \ra \infty}
\lambda_m = \lambda^{\bp}_{L,H}$ since the support of any smooth function with compact support on $H \setminus \Sigma$ is contained in $D_m$ for $m$ sufficiently large, cf.\ also \cite[Theorem 4.4.1]{P}.\\

Next we modify the operator $-L(u) = \sum_{i,j} a_{ij} \cdot \frac{\p^2 u}{\p x_i \p x_j} + \sum_i b_i \cdot \frac{\p u}{\p x_i} + c \cdot u$ by replacing $c$ with a new $C^\beta$-function $c[m]$ for a given $m$: We slightly decrease the value of $c$ within a small ball $B\subset D_1$ while keeping it fixed outside $B$, and such that the operator
\[
-L[m](u) = \sum_{i,j} a_{ij} \cdot \frac{\p^2 u}{\p x_i \p x_j} + \sum_i b_i \cdot \frac{\p u}{\p x_i} + c[m] \cdot u
\]
is still \si-adapted. Further, in view of the Rayleigh quotient and the fact that $\lambda_m  > \lambda^{\bp}_{L,H}$ we choose $c[m]$ such that
\begin{itemize}
  \item the first eigenvalue $\lambda_m[m]$ for the Dirichlet problem on $\overline{D_m}$ of the operator $\delta_{\bp}^2 \cdot L[m]$ just becomes $\lambda_m[m]=\lambda^{\bp}_{L,H}$.
  \item $c[m] \ra c$ in $C^{\beta}$-norm as $m \ra \infty$.
\end{itemize}
The resulting sequence of first eigenfunctions $\phi_m[m]$ with $\phi_m[m](p_0) = 1$ contains a compactly converging subsequence with limit $\phi >0$ on $H \setminus\Sigma$. This limit solves the equation $L  \, f = \lambda^{\bp}_{L,H} \cdot \bp^2 \cdot f$ on $H \setminus \Sigma$. Indeed, since $c[m] \ra c$ as $m \ra \infty$, we get Harnack inequalities for the $\phi_m[m]$ on domains $D \subset \overline{D} \subset H \setminus \Sigma$ for sufficiently large $m$ with constants independent of $m$, cf.\ \cite[Theorem 8.20]{GT} or \cite[Theorem 4.4.1]{P}.\\

It remains to show that $\phi$ is the unique positive solution up to multiples. So assume that $\psi>0$ is another solution of $L \, f = \lambda^{\bp}_{L,H} \cdot\bp^2 \cdot f$ on $H \setminus \Sigma$. Then we can find a constant $k>0$ such that $k \cdot \psi \ge \phi_m[m]$ near $\overline{B}$ for any $m$. But then $k \cdot \psi \ge\phi_m[m]$ on the whole of $H \setminus \Sigma$. Otherwise we can find, for any given $m$, a largest $k_0$ such that $k_0\cdot \psi \ge\phi_m[m]$, and there is some point $p \in  D_m \setminus \overline{B}$ with $k_0 \cdot \psi(p) = \phi_m[m](p)$. Hence, $k_0 \cdot \psi -\phi_m[m]$ would be a non-negative solution of $L  \, f = \lambda^{\bp}_{L,H} \cdot \bp^2 \cdot v$ on $D_m \setminus \overline{B}$ with an interior zero and positively lower bounded boundary values. However, this contradicts the \emph{Hopf maximum principle}, cf.\  \cite[Chapter 3.2]{GT}. (Note that in this case the maximum principle applies without the typical assumptions on the sign of the zeroth order coefficient as mentioned in the discussion after the proof of \cite[Theorem 3.5]{GT}. Actually, this follows from the strict inequality for the outer normal derivative $\p (k_0 \cdot \psi -\phi_m[m])/\p n> 0$ \cite[p.\ 24]{GT}, e.g., in extremal points of the zero set.) Hence $k_0 \cdot \psi -\phi_m[m] \equiv 0$ on the path component of $p$. It follows that $k \cdot \psi \ge \phi$ for some suitable $k>0$. Now we choose the smallest such $k$, which we call $k_1$. Then, either $k_1 \cdot \psi \equiv \phi$ or $k_1 \cdot
\psi > \phi$. In the latter case, we write $F:=k_1 \cdot \psi - \phi >0$ and repeat the previous argument to find some $l>0$ with $l\cdot F \ge \phi$. But then we also have $k_1\cdot \frac{l}{l+1} \cdot\psi \ge \phi$ contradicting the definition of $k_1$. Thus the solution $\phi$ is uniquely determined up to positive multiples.\\

\textbf{Supercritical Case} ($\lambda > \lambda^{\bp}_{L,H}$)  \,
In this case $L\, f = \lambda \cdot \bp^2 \cdot f$ has no positive solution. Indeed, assume to the contrary we had such a solution $u >0$. Then, we could find a smoothly bounded domain $D \subset \overline{D} \subset H \setminus \Sigma$ such that the eigenvalue $\lambda_D$ for the first Dirichlet eigenfunction $\phi_D$ of the operator $\delta_{\bp}^2 \cdot L$ on $D$ equals $\lambda$. Hence, there would exist a constant $k>0$ with $k \cdot u \ge \phi_D$. For the smallest such $k$ we would get an interior zero and positively lower bounded boundary values. As before this contradicts the Hopf maximum principle.\qed

%
%
%
\footnotesize
\renewcommand{\refname}{\fontsize{14}{0}\selectfont 
\end{document}